\documentclass{amsart}

\usepackage[foot]{amsaddr}
\let\oldtocsection=\tocsection

\let\oldtocsubsection=\tocsubsection 

\let\oldtocsubsubsection=\tocsubsubsection

\renewcommand{\tocsection}[2]{\vspace{0.5em}\hspace{0em}\oldtocsection{#1}{#2}}
\renewcommand{\tocsubsection}[2]{\vspace{0.5em}\hspace{1em}\oldtocsubsection{#1}{#2}}
\renewcommand{\tocsubsubsection}[2]{\vspace{0.5em}\hspace{2em}\oldtocsubsubsection{#1}{#2}}
\usepackage{graphicx,cancel,xcolor,hyperref,comment,graphicx,geometry}

\usepackage{tikz}
\usetikzlibrary{decorations.pathreplacing}
\usepackage{graphicx,url,etoolbox}
\usepackage{lipsum}
\usepackage{dsfont}
\usepackage{amsmath}
\usepackage{bbm}

\usepackage{verbatim}

\usepackage{fancyhdr}
\usepackage{lastpage}

\newtheorem{theoreme}{Theorem}[section]

\theoremstyle{definition}

\numberwithin{equation}{section}

\renewenvironment{proof}{{\bfseries \noindent Proof.}}{\demo}
\newcommand\xqed[1]{%
	\leavevmode\unskip\penalty9999 \hbox{}\nobreak\hfill
	\quad\hbox{#1}}
\newcommand\demo{\xqed{$\square$}}

\hypersetup{bookmarks, colorlinks, urlcolor=blue, citecolor=blue, linkcolor=red, hyperfigures, pagebackref,
	pdfcreator=LaTeX, breaklinks=true, pdfpagelayout=SinglePage, bookmarksopen=true,bookmarksopenlevel=2}

\usepackage{comment}

\def\R{\mathbb R}

\def\N{\mathbb N}

\def\HH{\mathcal H}

\def\AA{\mathcal A}

\def\la {{\lambda}}

\newcommand {\nc}   {\newcommand}
\nc {\be}   {\begin{equation}} \nc {\ee}   {\end{equation}} \nc
{\beq}  {\begin{eqnarray}} \nc {\eeq}  {\end{eqnarray}} \nc {\beqs}
{\begin{eqnarray*}} \nc {\eeqs} {\end{eqnarray*}}
\def\edc{\end{document}}

\providecommand{\abs}[1]{\lvert#1\rvert}

\numberwithin{equation}{section}

\theoremstyle{Thm}
\newtheorem{Thm}{Theorem}[section]
\newtheorem{lem}{Lemma}[section]
\newtheorem{prop}{Proposition}[section]

\definecolor{carnelian}{rgb}{0.7, 0.11, 0.11}
\definecolor{carmine}{rgb}{0.59, 0.0, 0.09}
\definecolor{burgundy}{rgb}{0.5, 0.0, 0.13}
\definecolor{darkmidnightblue}{rgb}{0.0, 0.2, 0.4}
\definecolor{dimgray}{rgb}{0.75, 0.75, 0.75}
\definecolor{palecarmine}{rgb}{0.69, 0.25, 0.21}
\newcounter{dummy} 
\numberwithin{dummy}{section}
\newtheorem{Theorem}[dummy]{Theorem}

\newtheorem{defi}[dummy]{Definition}

\numberwithin{equation}{section}

\def\AA{\mathcal A}
\def\HH{\mathbf{\mathcal H}}

\newcommand{\h}{\mathsf{h}}

\providecommand{\abs}[1]{\lvert#1\rvert}
\begin{document}
	\title[\fontsize{7}{9}\selectfont  ]{Polynomial stability of Thermoelastic Timoshenko system with non-global time-delayed Cattaneo's law}
\author{Haidar Badawi$^{1,2}$}
\author{Hawraa Alsayed$^{1,2}$}
\address{$^1$ PDEs with applications in materials sciences and biology, Department of Mathematics and Physics, Lebanese International University LIU, Beirut, Lebanon.}
\address{$^2$ PDEs with applications in materials sciences and biology, Department of Mathematics and Physics, Lebanese International University of Beirut BIU, Beirut, Lebanon.}
\email{haidar.badawi@liu.edu.lb, hawraa.sayed@liu.edu.lb}
\begin{abstract}
In this paper, we consider a one dimensional thermoelastic Timoshenko system in which the heat flux is given by Cattaneo's law and acts locally on the bending moment with a time delay. We prove its well-posedness, strong stability, and polynomial stability.
\end{abstract}
\maketitle
\pagenumbering{roman}
\maketitle
\pagenumbering{arabic}
\setcounter{page}{1}
\section{Introduction}
\noindent In this paper, we investigate the stability of a thermoelastic Timoshenko system in which the heat flux given by Cattaneo's law acts locally on the bending moment with a time delay. More precisely, we consider the following system:
\begin{equation}\label{pm-sysorig}
\left\{	\begin{array}{llll}
	\displaystyle \rho_1 \varphi_{tt}-k_1 (\varphi_x+\psi)_x =0 &\text{in} \ \  (0,\ell) \times (0,\infty),&\vspace{0.15cm}\\
	\displaystyle \rho_2 \psi_{tt}-k_2 \psi_{xx} +k_1 (\varphi_x +\psi)+d(\theta)=0&\text{in} \ \   (0,\ell) \times (0,\infty),\vspace{0.15cm}\\
	\rho_3 \theta_t +q_x +\delta \psi_{tx}=0 & \text{in} \ \ (0,\ell_0)\times (0,
	\infty), \, \ell_0 \in (0,\ell),\vspace{0.15cm}\\
	\gamma q_t +\mu_1 q+\mu_2 q(x,t-\tau)+\theta_x=0& \text{in} \ \ (0,\ell_0 )\times (0,\infty),\\ 
	\end{array}
	\right.
	\end{equation}
with boundary and initial conditions
\begin{equation}\label{pm-initialcon'}
	\left\{\begin{array}{llll}
		\varphi(0,t)=\varphi(\ell,t)=\psi(0,t)=\psi(\ell,t)=q(0,t)=\theta(\ell_0 ,t)=0 &\text{in } (0,\infty),&\vspace{0.15cm}\\	\varphi(x,0)=\varphi_0 (x),\ \varphi_t (x,0)=\varphi_1 (x), \ \psi(x,0)=\psi_0 (x)&\text{in } (0,\ell),& \vspace{0.15cm}\\
		\psi_t (x,0)=\psi_1 (x),\, \theta(x,0)=\theta_0 (x),\  q (x,0)=q_0 (x)&\text{in } (0,\ell_0),& \vspace{0.15cm}\\
		q(x,t) =f_0 (x,t)& \text{in } (0,\ell_0)\times (-\tau,0),& 
	\end{array}\right.\end{equation}
where $\varphi$ is the transverse displacement of the beam, $\psi$ is the rotational angle of a filament, $\theta$ is the temperature deviations, $q$ is the heat flux. The coefficients
 $\rho_1$, $\rho_2$, $\rho_3$, $k_1$, $k_2$,  $\gamma$, $\mu_1$, $\delta$ are  positive real numbers, $\mu_2$ is a non-zero real number. We define 
\begin{equation}\label{p3-a}
	d(\theta)=\left\{\begin{array}{lll}
\delta\theta_x& \text{in} &  (0,\ell_0 )\times(0,\infty),\vspace{0.15cm}\\
		0 & \text{in} &  (\ell_0,\ell )\times (0,\infty).
	\end{array}
	\right.
\end{equation}\label{sysintro}
In \cite{Timoshenko1921LXVIOT}, the classical Timoshenko-beam system is expressed by
\begin{equation*}
\left\{	\begin{array}{lll}
		\rho_1 \varphi_{tt}-k_1 (\varphi_x+\psi)_x=0 & \text{in} \  (0,\ell)\times(0,\infty),&\vspace{0.25cm}\\
		\rho_2 \psi_{tt}-k_2 \psi_{xx}+k_1 (\varphi_x+\psi)=0& \text{in} \  (0,\ell)\times(0,\infty).&
	\end{array}\right.
\end{equation*}
There are several publications concerning the stabilization of Timoshenko system  with different kinds of damping (see  \cite{Akil2020}, \cite{BASSAM20151177}, \cite{doi:10.1002/zamm.201500172} and  \cite{doi:10.1080/00036810903156149}).
Alves and al.  \cite{Alves2020} studied the thermoelastic Timoshenko system by coupling the thermal laws on both the shear force and the bending moment under Fourier's law, by considering:
\begin{equation}\label{sysintro'}
	\left\{	\begin{array}{llll}
		\displaystyle \rho_1 \varphi_{tt}-k_1 (\varphi_x+\psi)_x+\delta_1 \theta_x =0 &\text{in} \ \  (0,\ell) \times (0,\infty),&\vspace{0.15cm}\\
		\displaystyle \rho_2 \psi_{tt}-k_2 \psi_{xx} +k_1 (\varphi_x +\psi)-\delta_1 \theta+\delta_2 \vartheta_x=0&\text{in} \ \   (0,\ell) \times (0,\infty),\vspace{0.15cm}\\
		\rho_3 \theta_t -\mu_1 \theta_{xx} +\delta_1 (\varphi_x+\psi)_t=0& \text{in} \ \ (0,\ell)\times (0,
		\infty),\vspace{0.15cm}\\
		\rho_4 \vartheta_t -\mu_2 \vartheta_{xx} +\delta_2 \psi_{tx}=0 & \text{in} \ \ (0,\ell)\times (0,
		\infty),\vspace{0.15cm}\\
	\end{array}
	\right.
\end{equation}
they proved an exponential stability result without any relation between the coefficients of the system. Moreover,
\begin{itemize}
\item Fernandez Sare and Racke \cite{FernandezSare2009} studied system \eqref{sysintro'} by taking $\delta_1 =0$ and they showed that the system is exponentially stable if and only if \begin{equation}\label{speed}
	\frac{k_1}{\rho_1}=\frac{k_2}{\rho_2}
	.\end{equation}\\
\item  Almeida Junior and al. \cite{AlmeidaJunior2014}  studied system \eqref{sysintro'} by taking $\delta_2=0$ and they also proved exponential stability result provided that \eqref{speed} holds.
\end{itemize}
 In \cite{Djellali2021}, the authors studied the thermoelastic Timoshenko system where the heat flux given by  Cattaneo's law is acting on both  shear force and the bending moment, by considering:
\begin{equation}\label{sysintro''}
	\left\{	\begin{array}{llll}
		\displaystyle \rho_1 \varphi_{tt}-k_1 (\varphi_x+\psi)_x+\delta_1 \theta_x =0 &\text{in} \ \  (0,\ell) \times (0,\infty),&\vspace{0.15cm}\\
		\displaystyle \rho_2 \psi_{tt}-k_2 \psi_{xx} +k_1 (\varphi_x +\psi)-\delta_1 \theta+\delta_2 \vartheta_x=0&\text{in} \ \   (0,\ell) \times (0,\infty),\vspace{0.15cm}\\
		\rho_3 \theta_t +q_x +\delta_1 (\varphi_x+\psi)_t=0 & \text{in} \ \ (0,\ell)\times (0,
		\infty),\vspace{0.15cm}\\
		\gamma q_t + q+\mu_1\theta_x=0& \text{in} \ \ (0,\ell)\times (0,\infty),\vspace{0.15cm}\\ 
			\rho_4 \vartheta_t +p_x +\delta_2 \psi_{tx}=0 & \text{in} \ \ (0,\ell)\times (0,
		\infty),\vspace{0.15cm}\\
		\gamma_2 p_t + p+\mu_2\vartheta_x=0 & \text{in} \ \ (0,\ell)\times (0,\infty),\\ 
	\end{array}
	\right.
\end{equation}
with boundary conditions
$$
\varphi_x(0,t)=\varphi_x (\ell,t)=\psi(0,t)=\psi(\ell,t)=\theta(0,t)=\theta(\ell,t)=p(0,t)=p(\ell,t) \ \text{in} \ (0,\infty),
$$
they established an exponential stability result irrespective of the values of the coefficients of the system. However, Fernandez Sare and Racke \cite{FernandezSare2009} considered system \eqref{sysintro''} with $\delta_1=0$ and they showed that the system is not exponentially stable even if \eqref{speed} holds.\\\linebreak 
The originality of this work lies in the study of the thermoelastic Timoshenko system \eqref{pm-sysorig}-\eqref{pm-initialcon'}, in which the heat flux given by Cattaneo's law acts only on the bending moment and is distributed locally  with a time delay.\\\linebreak
In fact, from a physical point of view, models with partial differential equations that take into account a certain time delay are generally more realistic than those that do not take into account the effect of the delay, we refer the reader to  \cite{Akil2022, badawicpaa} for systems with time delay.\\\linebreak
This paper is organized as follows. In Section \ref{p3-WPS}, we prove the well-posedness of our system by using semigroup approach. In Section \ref{p3-sec3}, following a general criteria of Arendt and Batty, we show the strong stability of our system in the absence of the compactness of the resolvent. Finally, in Section \ref{p3-secpoly}, by combining the frequency domain approach with a specific multiplier multiplier method, we prove that the energy of our system decays polynomially with the rate $t^{-1}$.

\section{Well-posedness of the system}\label{p3-WPS}
\noindent In this section, we will establish the well-posedness of  system \eqref{pm-sysorig}-\eqref{pm-initialcon'} by using semigroup approach. For this aim, as in \cite{Nicaise2006},
we introduce the following auxiliary change of variable
\begin{equation}
z(x,s,t):=q (x,t-s\tau) \quad \text{in} \ \ (0,\ell_0 )\times(0,1)\times (0,\infty).
\end{equation}
Then, system \eqref{pm-sysorig}-\eqref{pm-initialcon'} becomes
\begin{equation}\label{pm-sysorig1}
	\left\{	\begin{array}{llll}
		\displaystyle \rho_1 \varphi_{tt}-k_1 (\varphi_x+\psi)_x =0 &\text{in} \ \  (0,\ell) \times (0,\infty),&\vspace{0.15cm}\\
		\displaystyle \rho_2 \psi_{tt}-k_2 \psi_{xx} +k_1 (\varphi_x +\psi)+d(\theta)=0&\text{in} \ \   (0,\ell) \times (0,\infty),\vspace{0.15cm}\\
		\rho_3 \theta_t +q_x +\delta \psi_{tx}=0 & \text{in} \ \ (0,\ell_0)\times (0,
		\infty),\vspace{0.15cm}\\
		\gamma q_t +\mu_1 q+\mu_2 z(x,1,t)+\theta_x=0 & \text{in} \ \ (0,\ell_0)\times (0,\infty),\vspace{0.15cm}\\ 
		\tau z_t (x,s,t)+z_s (x,s,t)=0& \text{in} \ \ (0,\ell_0)\times (0,1)\times (0,\infty),
	\end{array}
	\right.
\end{equation}
with the following boundary conditions 
\begin{equation}\label{pm-2.3}
\varphi(0,t)=\varphi(\ell,t)=\psi(0,t)=\psi(\ell,t)=q(0,t)=\theta(\ell_0 ,t)=0 \ \ \text{in} \ \ (0,\infty)
\end{equation}
and the following initial conditions
\begin{equation}\label{pm-initialcon}
\left\{\begin{array}{llll}
\varphi(x,0)=\varphi_0 (x),\ \varphi_t (x,0)=\varphi_1 (x), \ \psi(x,0)=\psi_0 (x)&\text{in } (0,\ell),& \vspace{0.15cm}\\
 \psi_t (x,0)=\psi_1 (x),\, \theta(x,0)=\theta_0 (x),\  q (x,0)=q_0 (x)&\text{in } (0,\ell_0),& \vspace{0.15cm}\\
z(x,s,0) =f_0 (x,-s \tau)& \text{in } (0,\ell_0)\times (0,1).& 
\end{array}\right.\end{equation} 
We set
$$
\Phi=(\varphi,\varphi_t,\psi,\psi_t,\theta,q,z)^\top,
$$
then system \eqref{pm-sysorig1}-\eqref{pm-initialcon} can be written as the following first order evolution equation
\begin{equation}\label{pm-Ut}
\Phi_t=\AA \Phi, \ \Phi(0)=\Phi_0=(\varphi_0,\varphi_1,\psi_0,\psi_1,\theta_0,q_0,f_0(\cdot,-\cdot \tau))^\top \in \HH,
\end{equation}
where the Hilbert space $\HH$ defined by 
$$
\HH=[H^1_0 (0,\ell)\times L^2 (0,\ell)]^2\times [L^2 (0,\ell_0)]^2 \times L^2 ((0,\ell_0 )\times (0,1)),
$$
is endowed with the inner product 
$$
\begin{array}{lll}
(\Phi,\Phi^1)_\HH&=&\displaystyle \int_0^\ell  \left(k_1 (\varphi_x+\psi)\overline{(\varphi_x^1+\psi^1)}+\rho_1 u\overline{u^1}+k_2\psi_x \overline{\psi_x^1}+\rho_2 v\overline{v^1}\right)dx+\int_{0}^{\ell_0}(\rho_3 \theta \overline{\theta^1}+\gamma q\overline{q})dx\vspace{0.15cm}\\
&&\displaystyle +\, \tau |\mu_2 |\int_{0}^{\ell_0}  \int_0^1 z \overline{z^1 }ds dx,
\end{array}
$$
where $\Phi=(\varphi,u,\psi,v,\theta,q,z)^\top, \Phi^1=(\varphi^1,u^1,\psi^1,v^1,\theta^1,q^1,z^1)^\top \in \HH$. Besides, 
the operator $\AA: D(\AA) \subset \HH \longmapsto \HH$ is defined by 
\begin{equation}\label{p2-op}
	\AA\begin{pmatrix}
		\varphi\\u\\\psi\\v\\\theta\\ q\\ z
	\end{pmatrix}=
	\begin{pmatrix} 
		\displaystyle u\vspace{0.15cm}\\
		\displaystyle k_1 \rho_1^{-1}(\varphi_x+\psi)_x\vspace{0.15cm}\\
		\displaystyle v\vspace{0.15cm}\\
		\displaystyle \rho_2^{-1}[ k_2\psi_{xx}-k_1(\varphi_x+\psi)- d(\theta)]\vspace{0.15cm}\\
		-\rho_3^{-1}(q_x +\delta v_x)\vspace{0.15cm}\\
		\displaystyle -\gamma^{-1} [\mu_1q +\mu_2 z(\cdot,1)+\theta_x]\vspace{0.15cm}\\
		\displaystyle -\tau^{-1}z_s
	\end{pmatrix},
\end{equation}
with domain
\begin{equation}
	D(\AA)=\left\{\begin{array}{lll}\vspace{0.25cm}
		\Phi=(\varphi,u,\psi,v,\theta,q,z)^{\top}\in \HH \ \  \text{such that} \\\vspace{0.25cm}\varphi, \psi  \in H^2 (0,\ell), \  u, v \in H_0^1 (0,\ell), \\\vspace{0.25cm}
		\displaystyle     \theta \in H^1_R (0,\ell_0):=\left\{\theta \in H^1 (0,\ell_0)\ | \ \theta (\ell_0)=0 \right\},\\  q \in H^1_L (0,\ell_0):=\left\{q \in H^1 (0,\ell_0)\ | \ q (0)=0 \right\}, \vspace{0.25cm}\\ z  \in L^2((0,\ell_0),H^1(0,1)), \  z(\cdot,0)=q \ \text{in} \ (0,\ell_0)
	\end{array}\right\} 
\end{equation}
The energy of system \eqref{pm-sysorig1}-\eqref{pm-initialcon} is given by 
$$
E(t)=\frac{1}{2}\left \|\left(\varphi,\varphi_t,\psi,\psi_t,\theta,q,z\right)^\top \right\|_{\HH}^2.
$$
Here and throughout the paper, the assumptions on $\mu_1$ and $\mu_2$ satisfy
\begin{equation*}
 \abs{\mu_2}<\mu_1.
\end{equation*}

\begin{prop}\label{p3-mdissip}
	{\rm  The operator $\AA$ defined above is m-dissipative.}
\end{prop}
\begin{proof}
For all $\Phi=(\varphi,u,\psi,v,\theta,q,z)^\top \in D(\AA)$, we have
\begin{equation*}\label{p3-Reau}
	\begin{array}{lll}
\displaystyle	\Re (\AA \Phi,\Phi)_{\HH}&=& \displaystyle \Re \left[ k_1\int_{0}^{\ell}  (u_x+v)(\overline{\varphi_x}+\overline{\psi})dx+k_1\int_{0}^{\ell} (\varphi_x+\psi)_x \overline{u}dx+k_2 \int_{0}^{\ell} v_x \overline{\psi}_xdx  +k_2 \int_{0}^{\ell} \psi_{xx}\overline{v}dx\vspace{0.15cm}\right.\vspace{0.15cm}\\\displaystyle &&-\,\displaystyle  k_1\int_{0}^{\ell}(\varphi_x +\psi)\overline{v}dx-\delta \int_{0}^{\ell_0 }\theta_x \overline{v}dx -\int_{0}^{\ell_0}(q_x+\delta v_x)\overline{\theta}dx-\mu_1 \int_{0}^{\ell_0}|q|^2 dx -\mu_2 \int_{0}^{\ell_0}z(\cdot,1)\overline{q}dx  \vspace{0.15cm} \vspace{0.15cm}\\\displaystyle 
&&\left. \displaystyle-\, \int_{0}^{\ell_0}\theta_x \overline{q}dx- \displaystyle  |\mu_2|\int_{0}^{\ell_0}|z(\cdot,1)|^2dx+|\mu_2|\int_{0}^{\ell_0}|q|^2 dx \right],
	\end{array}
\end{equation*}
using integration by parts and the fact that $\Phi \in D(\AA)$, we get 
\begin{equation*}
\Re (\AA \Phi,\Phi)_{\HH}= -\mu_1 \int_{0}^{\ell_0}|q|^2 dx -\Re \left[\mu_2 \int_{0}^{\ell_0}z(\cdot,1)\overline{q}dx \right] -\frac{|\mu_2|}{2}\int_{0}^{\ell_0}|z(\cdot,1)|^2dx+\frac{|\mu_2|}{2}\int_{0}^{\ell_0}|q|^2 dx.
\end{equation*}
Consequently, by Young's inequality, we claim 
\begin{equation}\label{pm-dissip}
\Re (\AA \Phi,\Phi)_\HH \leq -(\mu_1 -|\mu_2|)\int_{0}^{\ell_0 }|q|^2 dx\leq 0, 
\end{equation}
which implies that $\AA$ is dissipative. Let us prove that $\AA$ is maximal. For this aim, let $F=(f^1,f^2,f^3,f^4,f^5,f^6,f^7)^{\top}\in\HH$, we look for $\Phi=(\varphi,u,\psi,v,\theta,q,z)^{\top}\in D(\AA)$ a unique solution of 
\begin{equation}\label{pm-AU=F}
	-\AA \Phi=F\in \HH.
\end{equation}
Detailing \eqref{pm-AU=F}, we obtain 
\begin{eqnarray}
u&=&-f^1\in H^1_0 (0,\ell),\label{pm-f1}\\
-k_1\left(\varphi_x+\psi \right)_x &=&\rho_1 f^2\in L^2(0,\ell)\label{pm-f2}\\
v&=&-f^3\in H^1_0(0,\ell),\label{pm-f3}\\
	-k_2\psi_{xx} +k_1(\varphi_x+\psi)+d(\theta)&=&\rho_2 f^4\in L^2(0,\ell),\label{pm-f4}\\
q_x+\delta v_x &=&\rho_3 f^5\in L^2(0,\ell_0),\label{pm-f5}\\
\mu_1 q+\mu_2 z(\cdot,1)+\theta_x&=&\gamma f^6\in L^2(0,\ell_0),\label{pm-f6}\\
z_s&=&\tau f^7\in L^2((0,\ell_0)\times (0,1)), \label{pm-f7}
\end{eqnarray}
with the following boundary conditions
\begin{equation}\label{pm-bc}
\varphi (0)=\varphi(\ell)=\psi(0)=\psi(\ell)=q(0)=\theta(\ell_0)=0 \ \ \text{and} \ \ z(\cdot,0)=q \ \text{in} \ (0,\ell_0).
\end{equation}
From \eqref{pm-f7}, we obtain $z\in L^2( (0,\ell_0);H^1(0,1))$.
Moreover, from \eqref{pm-bc}, we deduce that 
\begin{equation}\label{pm-z}
	z=\tau \int_0^s f^7(\cdot,\xi_1 ) d\xi_1   +q.
\end{equation}
Inserting \eqref{pm-f3} in \eqref{pm-f5}, we get 
$$
q_x = \delta f^3_x +\rho_3 f^5,
$$
consequently, from \eqref{pm-bc} and the fact that $F\in \HH$, we deduce that $q\in H^1_L (0,\ell_0)$ and
\begin{equation}\label{pm-q}
 q=\delta f^3 +\rho_3 \int_0^x f^5(\xi_2) d\xi_2.
\end{equation}
Now, from \eqref{pm-f6}, we deduce that $\theta \in H^1_R (0,\ell_0)$ and 
\begin{equation}
\theta_x =-(\mu_1+\mu_2) \left(\delta f^3 +\rho_3 \int_0^x f^5(\xi_2) d\xi_2   \right)-\mu_2 \tau \int_0^1 f^7 (\cdot,\xi_1)d\xi_1  +\gamma f^6,
\end{equation}
as a consequence, we get
$$
\begin{array}{lll}
\displaystyle\theta= \int_{x}^{\ell_0} \left[(\mu_1+\mu_2)\left(\delta  f^3(\xi_3)  +\rho_3  \int_0^{\xi_3 } f^5(\xi_2) d\xi_2 \right)   
 +\mu_2 \tau \int_0^1 f^7 (\xi_3,\xi_1)d\xi_1-\gamma  f^6(\xi_3)\right]d\xi_3.
  \end{array}
$$
	Let $(\phi^1 ,\phi^2) \in \left[H^{1}_0 (0,\ell)\right]^2$, multiplying  \eqref{pm-f2} and \eqref{pm-f4} by $\overline{\phi^1}$ and $\overline{\phi^2}$ respectively, integrating over $(0,\ell)$, then using formal integrations by parts, we obtain
	\begin{equation}\label{p3-vf}
	\mathcal{B}((\varphi,\psi),(\phi^1,\phi^2))=\mathcal{L}(\phi^1,\phi^2), \ \ \forall (\phi^1,\phi^2)\in \left[H^1_0 (0,\ell)\right]^2,
	\end{equation}
	where
$$
\begin{array}{lll}
\displaystyle 	\mathcal{B}((\varphi,\psi),(\phi^1,\phi^2))=\displaystyle \int_0^\ell \left[ k_1 (\varphi_x+\psi)(\overline{\phi_x^1}+\overline{\phi^2})+k_2 \psi_x \overline{\phi^2_x }\right] dx
\end{array}
$$
and
$$
\begin{array}{lll}
\displaystyle \mathcal{L}(\phi^1,\phi^2)&=&\displaystyle \int_{0}^{\ell} \left(\rho_1 f^2 \overline{\phi^1}+\rho_2 f^4 \overline{\phi^2}\right)dx\\&&  \displaystyle+\, \delta \int_{0}^{\ell_0} \left[(\mu_1+\mu_2) \left(\delta f^3 +\rho_3 \int_0^x f^5(\xi_2) d\xi_2  \right)+\mu_2 \tau \int_0^1 f^7 (\cdot,\xi_1)d\xi_1  -\gamma f^6\right]dx.
\end{array}
$$
It is easy to see that,  $\mathcal{B}$ is a sesquilinear, continuous and coercive form on $\left[ H^{1}_0 (0,\ell)\right]^2  \times \left[H^{1}_0 (0,\ell)\right]^2 $ and $\mathcal{L}$ is an antilinear and continuous form on $\left[ H^{1}_0 (0,\ell)\right]^2 $. Then, it follows by Lax-Milgram theorem that \eqref{p3-vf} admits a unique solution $(\varphi,\psi)\in \left[H^{1}_0 (0,\ell)\right]^2 $. By taking test-functions $(\phi^1,\phi^2)\in [\mathcal{D}(0,\ell)]^2 $, we see that (\eqref{pm-f2}, \eqref{pm-f4}) hold in the distributional sense, from which we deduce that  $(\varphi,\psi)\in  \left[H^2 (0,\ell)\cap H^{1}_0 (0,\ell)\right]^2  $ and hence,  $\Phi \in D(\AA) $ is a unique solution of \eqref{pm-AU=F}. Accordingly, $\mathcal{A}$ is an isomorphism, and since $\rho\left(\mathcal{A}\right)$ is open set of $\mathbb{C}$ (see Theorem 6.7 (Chapter III) in \cite{Kato01}),  we easily get $R(\lambda I -\mathcal{A}) = {\mathcal{H}}$ for a sufficiently small $\lambda>0 $. This, together with the dissipativeness of $\mathcal{A}$, imply that   $D\left(\mathcal{A}\right)$ is dense in ${\mathcal{H}}$   and that $\mathcal{A}$ is m-dissipative in ${\mathcal{H}}$ (see Theorems 4.5, 4.6 in  \cite{Pazy01}). 
\end{proof}\\\linebreak
According to Lumer-Phillips theorem (see \cite{Pazy01}), Proposition \ref{p3-mdissip} implies that the operator $\AA$ generates a $C_{0}$-semigroup of contractions $e^{t\AA}$ in $\HH$ which gives the well-posedness of \eqref{pm-Ut}. Then, we have the following result:
\begin{Thm}{\rm
		For all $\Phi_0 \in \HH$,  system \eqref{pm-Ut} admits a unique weak solution $\Phi(t)=e^{t\AA}\Phi_0  \in C^0 (\R_+ ,\HH).
		$ Moreover, if $\Phi_0 \in D(\AA)$, then the system \eqref{pm-Ut} admits a unique strong solution $\Phi(t)=e^{t\AA}\Phi_0 \in C^0 (\R_+ ,D(\AA))\cap C^1 (\R_+ ,\HH).$}
\end{Thm}
\section{Strong Stability}\label{p3-sec3}
\noindent In this section, we will prove the strong stability of  system \eqref{pm-sysorig1}-\eqref{pm-initialcon}. The main result of this section is the following theorem.
\begin{theoreme}\label{p3-strongthm2}
	{\rm	The $C_0-$semigroup of contraction $\left(e^{t\AA}\right)_{t\geq 0}$ is strongly stable in $\HH$; i.e., for all $\Phi_0\in \HH$, the solution of \eqref{pm-Ut} satisfies 
		$$
		\lim_{t\rightarrow \infty}\|e^{t\AA}\Phi_0\|_{\HH}=0.
		$$}
\end{theoreme}\noindent  \noindent According to Theorem \ref{App-Theorem-A.2}, to prove Theorem \ref{p3-strongthm2}, we need to prove that the operator $\AA$ has no pure imaginary eigenvalues and $\sigma(\AA)\cap i\R $ is countable. The proof of these results is not reduced to the analysis of the point spectrum of $\AA$ on the imaginary axis since its resolvent is not compact. Hence
the proof of Theorem \ref{p3-strongthm2} is divided into the following two Lemmas.
\begin{lem}\label{p3-ker}
	{\rm For all $\la \in \R$, we have  $i\la I-\AA$ is injective i.e.,
		$$\ker(i\la I-\AA)=\{0\}.$$

	}
\end{lem}
\begin{proof}
	From Proposition \ref{p3-mdissip}, we have $0\in \rho (\AA)$. We still need to show the result for $\la \in \R^{*}$. For this aim, suppose that $\la\neq0$ and let $\Phi=(\varphi,u,\psi,v,\theta,q,z)^{\top}\in D(\AA)$ such that 
	\begin{equation}\label{p2-AU=ilaU}
	\AA \Phi=i\la \Phi.
	\end{equation}Equivalently, we have the following system
	\begin{eqnarray}
	u&=&i\la \varphi\label{p3-f1ker},
	\\\displaystyle  k_1 (\varphi_x+\psi)_x&=&i\la \rho_1  u\label{p3-f2ker},
	\\\displaystyle v&=&i\la \psi \label{p3-f3ker},
	\\\displaystyle k_2 \psi_{xx}-k_1 (\varphi_x+\psi)-d(\theta)&=&i\la \rho_2 \la v \label{p3-f4ker},
	\\\displaystyle q_x +\delta v_x&=&-i\la \rho_3 \theta \label{p3-f5ker},\vspace{0.15cm}
	\\\displaystyle \mu_1 q+\mu_2 z(\cdot,1)+\theta_x &=&-i\la \gamma q
	\vspace{0.15cm}\label{p3-f6ker},
		\\\displaystyle z_s&=&-i\la \tau z 
	\vspace{0.15cm}\label{p3-f7ker}.
	\end{eqnarray} 
	From  \eqref{pm-dissip} and \eqref{p2-AU=ilaU}, we obtain 
	\begin{equation}\label{p2-dissi=0}
	0=\Re \left(i\la \Phi,\Phi\right)_\HH=\Re\left(\AA \Phi,\Phi\right)_{\HH}\leq -(\mu_1 -|\mu_2|)\int_{0}^{\ell_0}|q|^2 dx \leq 0.
	\end{equation}
	Thus, we have
	\begin{equation}\label{p3-2.43}
q=0 \ \ \text{in} \ \ (0,\ell_0),
	\end{equation}
consequently, from \eqref{p3-f7ker} and the fact that $z(\cdot,0)=q$ in $(0,\ell_0)$, we deduce that 
	\begin{equation}\label{p3-3.10}
z=qe^{-i\la \tau s}=0 \  \ \text{in} \ \ (0,\ell_0)\times (0,1).
	\end{equation}
Inserting \eqref{p3-2.43} and \eqref{p3-3.10} in \eqref{p3-f5ker}, we get
$$
\theta_x =0 \ \ \text{in} \ \ (0,\ell_0) \ \ \text{and} \ \ d(\theta)=0  \ \ \text{in} \ \ (0,\ell),
$$ 
since $\theta(\ell_0)=0$, we infer
\begin{equation}\label{3.12}
\theta=0 \  \ \text{in} \ \ (0,\ell_0).
\end{equation}
Now, \eqref{p3-f3ker} and \eqref{p3-f5ker} imply that 
\begin{equation}
v_x =\psi_x =0 \ \ \text{in} \ \ (0,\ell_0).
\end{equation}
Since $\psi(0)=v(0)=0$, then
\begin{equation}\label{3.13}
v=\psi=0 \ \ \text{in } \ \ (0,\ell_0),
\end{equation}
From \eqref{p3-f4ker}, \eqref{3.12}, \eqref{3.13}, we obtain
$$
\varphi_x =0 \ \text{in} \ (0,\ell_0),
$$
hence, using $\varphi(0)=0$ and \eqref{p3-f1ker}, we get
$$
u=\varphi=0 \ \text{in} \ (0,\ell_0).
$$
Inserting \eqref{p3-f1ker} and \eqref{p3-f3ker} in \eqref{p3-f2ker} and \eqref{p3-f4ker}, we obtain 
\begin{equation}
\left\{\begin{array}{lll}
\displaystyle \varphi_{xx}=-\psi_{x}-\frac{\la^2 \rho_1}{k_1}\varphi & \text{in} &(\ell_0,\ell),\vspace{0.25cm}\\
\displaystyle \psi_{xx}=\frac{k_1}{k_2}\varphi_{x}+\frac{k_1-\la^2 \rho_2}{k_2}\psi & \text{in} & (\ell_0,\ell).
\end{array}\right.
\end{equation}
Let $\Psi=(\varphi,\varphi_{x},\psi,\psi_{x})^\top$. From the above system and the regularity of $\varphi$ and $\psi$, we deduce that 
\begin{equation}
\left\{\begin{array}{lll}
\Psi_x=A_\la \Psi \ \ \text{in} \ \ (\ell_0,\ell), \vspace{0.25cm}\\
\Psi(\ell_0)=0,
\end{array}\right.
\end{equation}
where 
$$
A_\la =\begin{pmatrix}
0&1&0&0\\
-\frac{\la^2 \rho_1}{k_1}&0&0&-1\\
0&0&0&1\\
0&\frac{k_1}{k_2}&\frac{k_1-\la^2\rho_2}{k_2}&0
\end{pmatrix}.
$$
The solution of the above differential equation is given by
$$
\begin{array}{lll}
\Psi (x)=\exp (A_\la (x-\ell_0))\Psi (\ell_0)=0, \ \forall x \in  (\ell_0,\ell).
\end{array}
$$
From the above estimation, \eqref{p3-f1ker} and \eqref{p3-f3ker}, we deduce that
$$
u=\varphi=v=\psi=0 \ \ \text{in} \ \ (\ell_0,\ell).
$$
Finally, we conclude that
$$
\Phi=0.
$$
\end{proof}
\begin{lem}\label{p3-surj}
	{\rm For all $\la \in \R $, we have 
	$$R(i\la I-\AA )=\HH.$$
	
}
\end{lem}
\begin{proof}
		From Proposition \ref{p3-mdissip}, we have $0\in\rho(\AA)$. We still need to show the result for $\la \in \R^{*}$. For this aim, let $F=(f^1,f^2,f^3, f^4,f^5,f^6, f^7)^{\top}\in \HH$, we look for $\Phi=(\varphi,u,\psi,v,\theta,q,z)^{\top}\in D(\AA)$ a solution of \begin{equation}\label{p3-ilaU-AU=F}
	(	i\la I -\AA)\Phi =F\in \HH.
	\end{equation}
	Detailing \eqref{p3-ilaU-AU=F}, we obtain
\begin{eqnarray}
	i\la \varphi-u=f^1\in H^1_0(0,\ell),\label{pm-f1'}\\
i\rho_1\la u	-k_1\left(\varphi_x+\psi \right)_x =\rho_1 f^2\in L^2(0,\ell) ,\label{pm-f2'}\\
i\la \psi	-v=f^3\in H^1_0(0,\ell),\label{pm-f3'}\\
i\rho_2 \la v 	-k_2\psi_{xx} +k_1(\varphi_x+\psi)+d(\theta)=\rho_2 f^4\in L^2(0,\ell),\label{pm-f4'}\\
i \rho_3 \la \theta+	q_x+\delta v_x =\rho_3 f^5\in L^2(0,\ell_0),\label{pm-f5'}\\
i\gamma \la q+	\mu_1 q+\mu_2 z(\cdot,1)+\theta_x=\gamma f^6\in L^2(0,\ell_0),\label{pm-f6'}\\
i\tau \la z+ 	z_s=\tau f^7\in L^2((0,\ell_0)\times (0,1)), \label{pm-f7'}
\end{eqnarray}
with the following boundary condition
$$
\varphi(0)=\varphi(\ell)=\psi(0)=\psi(\ell)=q(0)=\theta(\ell_0)=0 \ \text{and} \ z(\cdot,0)=q \ \text{in} \ (0,\ell_0).
$$
From \eqref{pm-f7'} and because $z(\cdot,0)=q \ \text{in} \ (0,\ell_0)$, we deduce that
\begin{equation}\label{z}
z=qe^{-i\tau \la s}+\tau \int_0^s e^{i\la \tau (\xi-s)}f^7(\cdot,\xi)d\xi.
\end{equation}
Consequently, from \eqref{pm-f6'}, we deduce that
\begin{equation}\label{thx}
\theta_x =-(i\gamma \la +\mu_1 +\mu_2 e^{-i\tau \la })q-\mu_2 \tau \int_0^1 e^{i\la \tau (\xi-1)}f^7(\cdot,\xi)d\xi+\gamma f^6.
\end{equation}
Moreover, since $\theta (\ell_0)=0$, then we get 
\begin{equation}\label{theta}
\theta=(i\gamma \la +\mu_1+\mu_2 e^{-i\tau \la})\int_x^{\ell_0}q dy +\mu_2 \tau \int_{x}^{\ell_0}\int_0^1 e^{i\la\tau(\xi-1)}f^7(y,\xi)d\xi dy-\gamma \int_{x}^{\ell_0}f^6dy.
\end{equation}
Inserting \eqref{pm-f3'} and the above equation in \eqref{pm-f5'}, we get 
\begin{equation}\label{3.27}
\begin{array}{lll}
\displaystyle i\rho_3 \la (i\gamma \la+\mu_1 +\mu_2 e^{-i\tau \la})\int_{x}^{\ell_0}q dy +q_x +i\delta \la \psi_x = G,
\end{array}
\end{equation}
where
$$
G=-i\rho_3 \mu_2 \tau \la \int_{x}^{\ell_0}\int_0^1 e^{i\la \tau (\xi-1)}f^7 (y,\xi)d\xi dy  +i\rho_3 \gamma \la \int_x^{\ell_0}f^6 dy+\rho_3 f^5 +\delta f^3_x.
$$
Inserting \eqref{pm-f1'} in \eqref{pm-f2'}, \eqref{pm-f3'} and \eqref{thx} in \eqref{pm-f4'}, then using the fact that 
$$
-q=\frac{d}{dx}\int_x^{\ell_0} q dy,
$$
we obtain
\begin{equation}\label{3.25}
	-\rho_1 \la^2 \varphi -k_1 (\varphi_x +\psi)_x =\rho_1 f^2 +i\rho_1 \la f^1,
\end{equation}
\begin{equation}\label{3.26}
	-\rho_2 \la^2 \psi -k_2 \psi_{xx}+k_1 (\varphi_x +\psi)+d_1\left( \frac{d}{dx}\int_x^{\ell_0}qdy\right)=\rho_2 f^4+i\rho_2 \la f^3+d_2 (f^6,f^7),
\end{equation}
where
$$
d_1\left( X\right)=\left\{\begin{array}{lll}
	\displaystyle (i\gamma \la +\mu_1 +\mu_2 e^{-i\tau \la }) X\vspace{0.15cm
	} & \text{in}  &(0,\ell_0),\\
	0 &  \text{in}& (\ell_0,\ell),
\end{array}\right.
$$
with $X= \frac{d}{dx}\int_x^{\ell_0}qdy$ or $\int_x^{\ell_0}qdy$,
and
$$
d_2(f^{6}, f^7)=\left\{\begin{array}{lll}
	\displaystyle \mu_2 \tau \int_0^1 e^{i\la \tau (\xi-1)}f^7(\cdot,\xi)d\xi-\gamma f^6\vspace{0.15cm
	} & \text{in}  &(0,\ell_0),\\
	0 &  \text{in}& (\ell_0,\ell).
\end{array}\right.
$$
We set 
$$
\begin{array}{l}
\mathbb{H}:=[H^1_0 (0,\ell)]^2 \times H^1_R (0,\ell_0), \ \ V:=(\varphi,\psi,\int_x^{\ell_0}qdy), \ \text{and} \  \ W:= (\Phi,\Psi,\int_x^{\ell_0}Qdy)\in \mathbb{H}.
\end{array}
$$
 Multiplying \eqref{3.25}, \eqref{3.26}, and \eqref{3.27} by $\overline{\Phi}$, $\overline{\Psi}$, and $ -\int_x^{\ell_0}\overline{Q}dy$, respectively, then integrating over $(0,\ell)$ and using formal integration by parts:
\begin{equation}\label{3.30}
	\begin{array}{l}
(\mathcal{B}_1+\mathcal{B}_2)(V,W)=\mathcal{L}(W), \ \forall W\in \mathbb{H},
\end{array}
\end{equation}
where
\begin{equation}
	\begin{array}{lll}
\displaystyle\mathcal{B}_1(V,W)&=&\displaystyle 	k_1 \int_{0}^{\ell}(\varphi_x +\psi)(\overline{\Phi_x}+\overline{\Psi})dx  +k_2 \int_0^{\ell} \psi_x \overline{\Psi_x}dx+\rho_3 \gamma \la^2 \int_0^{\ell_0} \left(\int_x^{\ell_0}qdy\right)\left(\int_x^{\ell_0}\overline{Q}dy\right)dx\vspace{0.15cm}\\
\displaystyle && -\,\displaystyle i\rho_3 \mu_1 \la \int_0^{\ell_0}\left(\int_x^{\ell_0}qdy\right)\left(\int_x^{\ell_0}\overline{Q}dy\right)dx +\int_0^{\ell_0}q \overline{Q}dx,
\end{array}
\end{equation}
\begin{equation}\label{B2}
	\begin{array}{lll}
		\displaystyle\mathcal{B}_2(V,W)&=& \displaystyle 	-\rho_1 \la^2 \int_{0}^{\ell}\varphi\overline{\Phi}dx    -\rho_2 \la^2 \int_0^{\ell}\psi \overline{\Psi}dx-\int_0^{\ell} d_1 \left(\int_x^{\ell_0}qdy\right)\overline{\Psi_x} dx\vspace{0.15cm}\\
		\displaystyle &&-\, \displaystyle i\rho_3 \mu_2 \la e^{-i\tau \la}\int_0^{\ell_0}\left(\int_x^{\ell_0}qdy\right)\left(\int_x^{\ell_0}\overline{Q}dy\right)dx +i\delta \la \int_0^{\ell_0} \psi \overline{Q}dx
	\end{array}
\end{equation}
and
\begin{equation}
\begin{array}{l}
	\mathcal{L}(W)=\displaystyle \int_0^{\ell}(\rho_1 f^2+i\rho_1\la f^1 )\overline{\Phi}dx +\int_0^{\ell}[\rho_2 f^4+i\rho_2 \la f^3+d_2 (f^6,f^7)]\overline{\Psi}dx -\int_0^{\ell_0} G\left(\int_x^{\ell_0}\overline{Q}dy \right) dx . 
	\end{array}
\end{equation}
Let $\mathbb{H}^{\prime}$  be the dual space of $\mathbb{H}$, define the following operators 
\begin{equation}
	\begin{array}{lll}
		\mathbb{B}_j :\ \ \mathbb{H}&\longmapsto&\mathbb{H}^{\prime} \\
		\	\ \ \ \ \ \, \ V&\longmapsto& \mathbb{B}_j V
	\end{array},\ \ j\in \{1,2\},
\end{equation}such that 
\begin{equation}\label{p2-3.63}
	\begin{array}{lll}
	 	\left(\mathbb{B}_j V\right)\left(W\right)=\mathcal{B}_j \left(V,W\right),\ \ \forall W \in \mathbb{H}.
	\end{array}
\end{equation}We need to prove that the operator $\mathbb{B}_1+\mathbb{B}_2$ is an isomorphism. For this purpose, we divide the proof into two steps:\\\linebreak
\textbf{Step 1.} In this step, we prove that the operator $\mathbb{B}_2 $ is compact. Indeed, from \eqref{B2} and \eqref{p2-3.63} we have 
\begin{equation}
	\begin{array}{l}
	\left|\mathcal{B}_2 \left(V,W\right)\right|\leq C\left\| V\right\|_{\left[L^2 (0,\ell)\right]^2\times L^2 (0,\ell_0)}\left\| W\right\|_\mathbb{H},
	\end{array}
\end{equation} 
for some positive constant $C$, and consequently, using the compact embedding from $\mathbb{H}$ into $\left[L^2 (0,\ell)\right]^2\times L^2 (0,\ell_0) $, we deduce that $\mathbb{B}_2 $ is a compact operator.
\\\linebreak
The compactness property and the fact that $\mathbb{B}_1$ is an isomorphism imply that the operator $\mathbb{B}_1+\mathbb{B}_2$ is a Fredholm operator of index 
$$
0=\dim \ker(\mathbb{B}_1+\mathbb{B}_2)-\dim R(\mathbb{B}_1+\mathbb{B}_2)^{\perp}.
$$
Now, following Fredholm alternative, we simply need to prove that the operator $\mathbb{B}_1+\mathbb{B}_2$ is injective to obtain that it is an isomorphism.
\\\linebreak
\textbf{Step 2.} In this step, we want to prove that the operator $\mathbb{B}_1+\mathbb{B}_2$ is injective (i.e. $\ker(\mathbb{B}_1+\mathbb{B}_2)=\{0\}$). So, let $\left(\widehat{\varphi},\widehat{\psi},\int_x^{\ell_0}\widehat{q}dy\right)\in \ker(\mathbb{B}_1+\mathbb{B}_2)$ which gives 
\begin{equation*}
	\begin{array}{l}
(	\mathcal{B}_1+\mathcal{B}_2)\left(\left(\widehat{\varphi},\widehat{\psi},\int_x^{\ell_0}\widehat{q}dy\right),W\right)=0,\ \ \forall W\in \mathbb{H}.
\end{array}
\end{equation*}
Thus, we find that 
\begin{equation}\label{3.25hat}
\left\{	\begin{array}{lll}
	-\rho_1 \la^2 \widehat{\varphi} -k_1 (\widehat{\varphi}_x +\widehat{\psi})_x =0\\
	-\rho_2 \la^2 \widehat{\psi} -k_2 \widehat{\psi}_{xx}+k_1 (\widehat{\varphi}_x +\widehat{\psi})+d_1\left( \frac{d}{dx}\int_x^{\ell_0}\widehat{q}dy\right)=0,\\
		\displaystyle i\rho_3 \la (i\gamma \la+\mu_1 +\mu_2 e^{-i\tau \la})\int_{x}^{\ell_0}\widehat{q} dy +\widehat{q}_x +i\delta \la \widehat{\psi}_x = 0,\\
		\widehat{\varphi}(0)=\widehat{\varphi}(\ell)=\widehat{\psi}(0)=\widehat{\psi}(\ell)=\widehat{q}(0)=0.
		\end{array}\right.
\end{equation}
Therefore, the vector $\widehat{\Phi}$ defined by $$
\begin{array}{l}
	\widehat{\Phi}=\left(\widehat{\varphi},i\la\widehat{\varphi},\widehat{\psi},i\la \widehat{\psi}, (i\gamma \la+\mu_1 +\mu_2 e^{-i\tau \la}) \int_x^{\ell_0}\widehat{q}dy,\widehat{q},\widehat{q}e^{-i\tau\la s}\right)^{\top}
\end{array}
$$ belongs to $D(\AA )$  and satisfies $$i\la \widehat{\Phi}-\AA\widehat{\Phi}=0,$$and consequently $\widehat{\Phi}\in \ker(i\la I-\AA)$. Hence, Lemma \ref{p3-ker} yields $\widehat{\Phi}=0$ as a result $\widehat{\varphi}=\widehat{\psi}=0$ and $\ker(\mathbb{B}_1+\mathbb{B}_2)=\{0\}$.\\\linebreak
Steps 1 and 2 guarantee that the operator $\mathbb{B}_1+\mathbb{B}_2$ is an isomorphism. Furthermore, it is easy to see that the operator $\mathcal{L}$ is an antilinear and continuous form on $\mathbb{H}$. Therefore, \eqref{3.30} admits a unique solution $V\in \mathbb{H}$. In \eqref{3.30}, by taking the test functions $W \in [\mathcal{D}(0,\ell)]^2\times \mathcal{D}(0,\ell_0)$, we see that (\eqref{3.27}, \eqref{3.25}, \eqref{3.26}) holds in distributional sense, hence 
$$
q\in H^1_L(0,\ell_0) \ \ \text{and} \  \ \varphi, \psi \in H^2 (0,\ell)\cap H^1_0 (0,\ell).
$$
Moreover, from \eqref{z} and \eqref{theta}, we deduce that 
$$
z\in L^2 ((0,\ell_0); H^1 (0,1))\ \ \text{and} \ \ \theta \in H^1_R (0,\ell_0).
$$
Finally, we claim that 
$$
\Phi =\left(\varphi,i\la \varphi-f^1,\psi, i\la \psi-f^3, \theta, q, z\right)^\top \in D(\AA),
$$
 is the unique solution of \eqref{p3-surj}.
\end{proof}
\\\linebreak
\textbf{Proof of Theorem \ref{p3-strongthm2}.} From Lemma \ref{p3-ker},  the operator $\AA $ has no pure imaginary eigenvalues (i.e. $\sigma_p (\AA)\cap i\R=\emptyset$). 
Moreover, from Lemma \ref{p3-ker} and Lemma \ref{p3-surj},   $i\la I-\AA $ is bijective for all $\lambda\in \mathbb{R}$ and since $\AA$ is closed, we conclude, with the help of the closed graph theorem, that  $i\la I-\AA $ is an isomorphism for all $\lambda\in \mathbb{R}$, hence
$\sigma(\AA )\cap i\R=\emptyset$.
\xqed{$\square$}

\section{Polynomial Stability }\label{p3-secpoly}
\noindent In this section, we will prove the polynomial stability of  system \eqref{pm-sysorig1}-\eqref{pm-initialcon}. The main result of this section is the following theorem. 

\begin{Thm}\label{p5-polthm}
	{\rm  For all $\Phi_0\in D(\mathcal{A})$, there exists a constant $C>0$ independent of $\Phi_0$ such that the energy of system \eqref{pm-sysorig1}-\eqref{pm-initialcon} satisfies the following estimation 
		$$
		E(t)\leq \frac{C}{t}\|\Phi_0\|^2_{D(\AA)},\quad \forall\, t>0.
		$$}
\end{Thm}
\noindent According to Theorem \ref{bt}, to prove Theorem \ref{p5-polthm}, we   need to prove  the following two conditions
\begin{equation}\label{p5-4.1}
	i\R \subset \rho (\AA)
\end{equation}and
\begin{equation}\label{p5-Condition-2-pol}
	\limsup_{\la\in \R,\ \abs{\la}\rightarrow \infty}\frac{1}{\la^2 }\|(i\la I-\AA)^{-1}\|_{\mathcal{L}(\HH)}<\infty.
\end{equation}
As condition \eqref{p5-4.1} was checked in Section \ref{p3-sec3}, we only need to prove the second condition. Condition \eqref{p5-Condition-2-pol} is proved by a contradiction argument. For this purpose,
suppose that \eqref{p5-Condition-2-pol} is false, then there exists $\{(\la_n,\Phi_n:=(\varphi_n,u_n,\psi_n,v_n,\theta_n, q_n, z_n)^{\top})\}_{n\geq1}\subset \R^{\ast}_{+}\times D(\mathcal{A})$ with
\begin{equation}\label{p5-contra-pol2}
	\la_n \to\infty \  \text{ as } n\to \infty \quad \text{and}\quad \|\Phi_n\|_{\mathcal{H}}=1, \forall n\geq 1,
\end{equation}
such that  
\begin{equation}\label{p3-eq0ps}
	(\la_n )^2 (i\la_n I-\AA )\Phi_n =F_n:=(f^1_n, f^2_n, f^3_n, f^4_n, f^5_n, f^6_n, f^7_n)^{\top}  \to 0  \quad \text{in}\quad \HH \quad \text{ as } n\to \infty.
\end{equation} 
For simplicity, we drop the index $n$. Equivalently, from \eqref{p3-eq0ps}, we have

\begin{eqnarray}
	i\la \varphi-u=\la^{-2}f^1,\label{pm-f1p}\\
	i\rho_1\la u	-k_1\left(\varphi_x+\psi \right)_x =\rho_1\la^{-2} f^2 ,\label{pm-f2p}\\
	i\la \psi	-v=\la^{-2}f^3,\label{pm-f3p}\\
	i\rho_2 \la v 	-k_2\psi_{xx} +k_1(\varphi_x+\psi)+d(\theta)=\rho_2\la^{-2} f^4,\label{pm-f4p}\\
	i \rho_3 \la \theta+	q_x+\delta v_x =\rho_3 \la^{-2}f^5,\label{pm-f5p}\\
	i\gamma \la q+	\mu_1 q+\mu_2 z(\cdot,1)+\theta_x=\gamma \la^{-2} f^6,\label{pm-f6p}\\
	i\tau \la z+ 	z_s=\tau \la^{-2}f^7. \label{pm-f7p}
\end{eqnarray}
Inserting \eqref{pm-f1p} in \eqref{pm-f2p}, \eqref{pm-f3p} in \eqref{pm-f4p} and \eqref{pm-f5p}, we get 
\begin{eqnarray}
\rho_1 \la^2 \varphi+k_1 (\varphi_x+\psi)_x &=&-\rho_1 \la^{-2}f^2-i\rho_1 \la^{-1}f^1,\label{4.10p}\\
\rho_2 \la^2 \psi+k_2 \psi_{xx}-k_1 (\varphi_x+\psi)-d(\theta)&=&-\rho_2 \la^{-2}f^4-i\rho_2 \la^{-1}f^3,\label{4.11p}\\
i\rho_3 \la \theta +q_x +i\delta \la \psi_x &=&\rho_3 \la^{-2}f^5 +\delta \la^{-2}f^3_x\label{4.12p}.
\end{eqnarray}
Here we will check condition \eqref{p5-Condition-2-pol} by finding a contradiction with \eqref{p5-contra-pol2} i.e. by showing $\left\|\Phi\right\|_{\HH}=o(1)$.  For clarity, we divide the proof into three Lemmas. From the above system and the fact that $\|\Phi\|_\HH=1$ and $\|F\|_\HH=o(1)$, we infer that 
\begin{equation}\label{star}
\left\{	\begin{array}{lll}
\|\la \varphi\|_{L^2(0,\ell)}, \|\la \psi\|_{L^2(0,\ell)}, \|\varphi_x\|_{L^2(0,\ell)}, \|\psi_x\|_{L^2(0,\ell)}, \|\theta\|_{L^2(0,\ell_0)}=O(1),\vspace{0.15cm}\\ \|\varphi_{xx}\|_{L^2(0,\ell)}, \|\psi_{xx}\|_{L^2(0,\ell)}=O(\la).
\end{array}\right.
\end{equation}
Also, for all $\Omega \subseteq (0,\ell)$, we have
\begin{equation}\label{star2}
\|\psi_{xx}\|_{L^2(\Omega)}\leq O(\la)\|\la \psi\|_{L^2(\Omega)}+O(1).
\end{equation}

\begin{lem}\label{p3-1stlemps}
		{\rm The solution $\Phi=(\varphi,u,\psi,v,\theta,q,z)^{\top}\in D(\AA)$ of  \eqref{pm-f1p}-\eqref{pm-f7p} satisfies the following estimations 
			\begin{eqnarray*}
		\int_0^{\ell_0}|q|^2dx =o(\la^{-2}), \\
		 \int_0^{\ell_0}\int_0^1 |z|^2 ds dx=o(\la^{-2}), \\ \int_0^{\ell_0} |z(\cdot,1)|^2 dx =o(\la^{-2}),\\
		 \int_0^{\ell_0} |\theta_x|^2 dx =o(1),\\
		\int_0^{\ell_0}|\theta|^2 dx =o(\la^{-1}).
		\end{eqnarray*}
	}
\end{lem}
\begin{proof}
		First, taking the inner product of \eqref{p3-eq0ps} with $\Phi$ in $\HH$ and using \eqref{pm-dissip}, we get
	\begin{equation*}\label{p3-4.10}
	\displaystyle	\int_0^{\ell_0}|q|^2 dx \leq \la^{-2} (\mu_1-|\mu_2|)^{-1}\Re (F,\Phi)_\HH\leq \la^{-2}(\mu_1-|\mu_2|)^{-1}\|F\|_\HH \|\Phi\|_\HH=o(\la^{-2}).
	\end{equation*}
 From  \eqref{pm-f7p} and since $z(\cdot,0)=q\ \text{in} \ (0,\ell_0)$, we have
	\begin{equation*}\label{p3-4.11}
	z=qe^{-i\tau \la s}+\tau \la^{-2}\int_0^s e^{i\la \tau (\xi-s)}f^7(\cdot,y)dy \ \ \text{in} \ \ (0,\ell_0)\times(0,1),
	\end{equation*}
consequently, by using Cauchy-Schwarz's inequalities, we deduce that 

\begin{equation*}
	\begin{array}{lll}
		\displaystyle \int_0^{\ell_0}\int_0^1 |z|^2 dsdx  &\leq & \displaystyle 2\int_0^{\ell_0} |q|^2 dx +2 \tau^2\la^{-4}\int_0^{\ell_0}\int_0^1 \left( \int_{0}^s |f^7 (\cdot,y)| dy\right)^2 dsdx\vspace{0.25cm}\\
	&\leq &\displaystyle  2\int_0^{\ell_0} |q|^2 dx +2\tau^2\la^{-4}\int_0^{\ell_0}\int_0^1 s \int_0^s |f^7(\cdot,y)|^2 dy ds dx\vspace{0.25cm}\\
&\leq &	\displaystyle  2\int_0^{\ell_0} |q|^2 dx +2 \tau^2\la^{-4}\left(\int_0^1 s ds\right)\int_0^{\ell_0} \int_0^1 |f^7(\cdot,y)|^2 dy  dx \vspace{0.25cm}\\
&=&	\displaystyle   2\int_0^{\ell_0}|q|^2 dx+\tau^2\la^{-4} \int_0^{\ell_0}\int_{0}^1 |f^7 (\cdot,y)|^2 dy dx\vspace{0.15cm}\\
		&=&\displaystyle o(\la^{-2}),
	\end{array}
\end{equation*}
and 
$$
\begin{array}{lll}
	\displaystyle \int_0^{\ell_0} |z(\cdot,1)|^2 dx&=& \displaystyle \int_0^{\ell_0}\left| qe^{-i\tau \la }+\tau \la^{-2}\int_0^1 e^{i\la \tau (\xi-1)}f^7(\cdot,y)dy \right|^2 dx \vspace{0.15cm}\\
	& \leq & \displaystyle 2  \int_0^{\ell_0} |q|^2 dx +2\tau^2\la^{-4}\int_0^{\ell_0} \left(\int_0^1 |f^7(\cdot,y)|dy\right)^2dx\vspace{0.25cm}\\
	& \leq& \displaystyle  2\int_0^{\ell_0} |q|^2 dx +2\tau^2\la^{-4}\int_0^{\ell_0} \int_0^1 |f^7(\cdot,y)|^2 dy dx\vspace{0.15cm}\\
	&=&o(\la^{-2}).
\end{array}
$$
Now, from \eqref{pm-f5p}, we obtain
$$
\begin{array}{lll}
\displaystyle \int_0^{\ell_0}|\theta_x|^2 dx &=& \displaystyle  \int_0^{\ell_0}   \left| 	-i\gamma \la q-	\mu_1 q-\mu_2 z(\cdot,1)+\gamma \la^{-2} f^6\right|^2dx\vspace{0.15cm}\\
  &\leq& \displaystyle  (\gamma^2 \la^2 +\mu_1^2)\int_0^{\ell_0}|q|^2 dx +\mu_2^2 \int_0^{\ell_0} |z(\cdot,1)|^2 dx +\gamma^2 \la^{-4}\int_{0}^{\ell_0} |f^6|^2 dx\vspace{0.15cm}\\
  &=&o(1). 
\end{array}
$$
Finally, multiplying \eqref{4.11p} by $ (\rho_3 \la)^{-1}$$\overline{\theta}$, integrating over $(0,\ell_0)$, using integration by parts:
$$
\begin{array}{lll}
\displaystyle \int_{0}^{\ell_0}|\theta|^2 dx &=& \displaystyle \Im\left[(\rho_3 \la)^{-1}\int_0^{\ell_0} q \overline{\theta_x}dx -i\delta \int_{0}^{\ell_0}\psi\overline{\theta_x}dx +\la^{-2}\int_{0}^{\ell_0} (\rho_3f^5+\delta f^3_x )\overline{\theta}dx\right]\vspace{0.15cm}\\\displaystyle 
&\leq& \displaystyle  (\rho_3 \la)^{-1}\|q\|_{L^2 (0,\ell_0)}\|\theta_x\|_{L^2 (0,\ell_0)}+\delta \underbrace{\|\psi\|_{L^2 (0,\ell_0)}}_{=O(\la^{-1})}\|\theta_x\|_{L^2 (0,\ell_0)}+\la^{-2}(\rho_3+\delta)\|F\|_{\HH}\underbrace{\|\theta\|_{L^2 (0,\ell_0)}}_{=O(1)}\vspace{0.15cm}\\
& =&o(\la^{-1}).
\end{array}
$$
\end{proof}
\begin{lem}\label{4.2}
	{\rm The solution $\Phi=(\varphi,u,\psi,v,\theta,q,z)^{\top}\in D(\AA)$ of  \eqref{pm-f1p}-\eqref{pm-f7p} satisfies the following estimations 
	$$
	\| \psi_x\|^2_{L^2(\omega_4)}=o\left(\la^{-2}\right) \ \ \text{and} \ \ \|\la\psi\|^2_{L^2(\omega_5)}=o(\la^{-2}).
	$$}
\end{lem}
\begin{proof}
	To prove this Lemma, for some $j,k \in \N$, $1\leq j\leq 7$, and $0<\varepsilon<\frac{\ell_0}{14}$, we need to define the following intervals:
	$$
	\omega_{k}=(k\varepsilon,\ell_0-k\varepsilon), \ \forall \  j-1 \leq k \leq j
	$$
	and as well we need to fix functions $\h_j\in C^{1}\left(\overline{\omega_0}\right)$ such that $0\leq \h_j (x)\leq 1$, for all $x\in\overline{\omega_0}$ and 
		\begin{equation*}
			\h_j (x)= 	\left \{ \begin{array}{lll}
				1 &\text{if} \quad \,\,  x \in \overline{\omega_j},&\vspace{0.1cm}\\
				0 &\text{if } \quad x \in \overline{\omega_0}\backslash \omega_{j-1}.&
			\end{array}	\right. \qquad\qquad
		\end{equation*}
\noindent For simplicity, let us divide the proof into four steps:\\\linebreak
	\textbf{Step1. } In this step, we prove:
	\begin{equation}\label{4.20p}
		\begin{array}{lll}
			\displaystyle \| \psi_x\|^2_{L^2(\omega_{j})}&\leq& \displaystyle  \max\left(\frac{\rho_3}{\delta} \max_{x\in \overline{\omega_j}} |\mathsf{h}_j^\prime(x)|, \frac{\rho_3} {\delta}, \frac{1}{\delta}\max_{x\in \overline{\omega_j}} |\mathsf{h}_j^\prime(x)|,\frac{1}{\delta},1 \right)\left[ \left(\|\theta\|_{L^2(\omega_{j-1})}+\|\theta_x\|_{L^2(\omega_{j-1})}\right)\|\psi\|_{L^2(\omega_{j-1})}\right.\vspace{0.15cm}\\
			 &&+\, \displaystyle \left. \la^{-1}\left(\|\psi_x\|_{L^2(\omega_{j-1})}+\|\psi_{xx}\|_{L^2(\omega_{j-1})}\right)\|q\|_{L^2(\omega_{j-1})}+\la^{-3}\|F\|_{\HH}\|\psi_x\|_{L^2(\omega_{j-1})}\right].
	\end{array}\end{equation}
For this aim, Multiplying \eqref{4.12p} by $\delta^{-1}\la^{-1}\h_{j}\overline{\psi_x}$ and integrating over $\omega_{j-1}$, we get
	\begin{equation}\label{4.24'}
	\int_{\omega_{j-1}}\h_{j} |\psi_x|^2dx  =\Im\left[-i\frac{\rho_3}{\delta}\ \int_{\omega_{j-1}}\h_{j} \theta \overline{\psi_x}dx -\frac{1}{\delta \la }\int_{\omega_{j-1}}\h_{j}q_x\overline{\psi_x}dx+\la^{-3}\int_{\omega_{j-1}}\h_{j} (\rho_3 \delta^{-1}f^5+f^3_x)\overline{\psi_x}dx \right].
	\end{equation}
	Since $\h_{j}(x)=0, \ \forall  x\in \partial \omega_{j-1}$, then integration by parts and Cauchy-Schwarz's inequality yield,
	$$
	\left|\int_{\omega_{j-1}}\h_{j}\theta\overline{\psi_x}dx\right| =\left|\int_{\omega_{j-1}}(\h_{j}^\prime \theta +\h_{j}\theta_x)\overline{\psi}dx \right|\leq \left(\max_{x\in \overline{\omega_j}} |\mathsf{h}_j^\prime(x)|\|\theta\|_{L^2(\omega_{j-1})}+\|\theta_x\|_{L^2(\omega_{j-1})}\right)\|\psi\|_{L^2(\omega_{j-1})}
	$$
	and 
	$$
	\left|\int_{\omega_{j-1}}\h_{j}q_x \overline{\psi_x}dx\right| =\left|\int_{\omega_{j-1}}(\h_{j}^\prime \overline{\psi_x}+\h_{j}\overline{\psi_{xx}})qdx\right|\leq \left(\max_{x\in \overline{\omega_j}} |\mathsf{h}_j^\prime(x)|\|\psi_x\|_{L^2(\omega_{j-1})}+\|\psi_{xx}\|_{L^2(\omega_{j-1})}\right)\|q\|_{L^2(\omega_{j-1})}. 
	$$
	Moreover, by Cauchy-Schwarz's inequality, we get
	$$
\left|	\int_{\omega_{j-1}}\h_{j} (\rho_3 \delta^{-1}f^5+f^3_x)\overline{\psi_x}dx\right| \leq  \left(\rho_3 \delta^{-1}+1\right)\|F\|_{\HH}\|\psi_x\|_{L^2(\omega_{j-1})}.
	$$
	Using \eqref{4.24'}, the above inequalities, and the fact that $\|\psi_x \|^2_{L^2 (\omega_j)}\leq \|\sqrt{\mathsf{h}_j}\psi_x \|^2_{L^2 (\omega_{j-1})} $, we deduce \eqref{4.20p}.\linebreak\\
	\textbf{Step 2.} In this step, we prove:
		\begin{equation}\label{4.21p}
		\begin{array}{lll}
			\displaystyle\| \la \psi\|^2_{L^2(\omega_{j})}&\leq&  \displaystyle  \max\left(k_2\rho_2^{-1} ,k_2\rho_2^{-1} \max_{x\in \overline{\omega_j}} |\mathsf{h}_j^\prime(x)|, k_1\rho_2^{-1} \max_{x\in \overline{\omega_j}} |\mathsf{h}_j^\prime(x)|, k_1\rho_2^{-1} , \delta \rho_2^{-1},\rho_2  \right)\left[\|\psi_x\|^2_{L^2(\omega_{j-1})}\right.\vspace{0.15cm}\\
			&&+\,  \left.\displaystyle \left(\|\psi_x\|_{L^2(\omega_{j-1})}+\|\theta_x\|_{L^2(\omega_{j-1})}\right)\|\psi\|_{L^2(\omega_{j-1})}\right.\vspace{0.15cm}\\&&+\,  \left.\displaystyle \left(\|\psi\|_{L^2(\omega_{j-1})}+\|\psi_x\|_{L^2(\omega_{j-1})}\right)\|\varphi\|_{L^2(\omega_{j-1})}\right.\vspace{0.15cm} \\&&+\,
			\left.\displaystyle  \|\psi\|^2_{L^2(\omega_{j-1})}+(\la^{-2}+\la^{-1})\|F\|_{\HH}\|\psi\|_{L^2(\omega_{j-1})}\right].
		\end{array}
	\end{equation}
So, multiplying \eqref{4.11p} by $\rho_2^{-1}\h_{j}\overline{\psi}$ and integrating over $\omega_{j-1}$, we infer
\begin{equation}\label{4.25}
	\begin{array}{lll}
		\displaystyle \int_{\omega_{j-1}}\h_{j}|\la \psi|^2 dx =\Re \left[-k_2 \rho_2^{-1}\int_{\omega_{j-1}}\h_j\psi_{xx}\overline{\psi}dx+k_1\rho_2^{-1}\int_{\omega_{j-1}}\h_{j}(\varphi_x+\psi)\overline{\psi}dx+\delta\rho_2^{-1}\int_{\omega_{j-1}}\h_{j}\theta_x \overline{\psi}dx \right.\vspace{0.15cm}\\\left.\hspace{3.5cm}
		\displaystyle -\rho_2\int_{\omega_{j-1}}\h_{j}( \la^{-2}f^4+i \la^{-1}f^3)\overline{\psi}dx\right].
	\end{array}
	\end{equation}
	Using integrations by parts and Cauchy-Schwarz's inequality, we get 
	$$
\left|	\int_{\omega_{j-1}}\h_{j}\psi_{xx}\overline{\psi}dx \right|=\left|-\int_{\omega_{j-1}}\h_{j}|\psi_x|^2dx -\int_{\omega_{j-1}}\h_{j}^\prime \psi_x \overline{\psi}dx \right|\leq  \|\psi_x\|^2_{L^2(\omega_{j-1})}+\max_{x\in \overline{\omega_j}} |\mathsf{h}_j^\prime(x)|\|\psi_x\|_{L^2(\omega_{j-1})}\|\psi\|_{L^2(\omega_{j-1})}
	$$
	and 
	$$
	\begin{array}{lll}
\displaystyle \left|	\int_{\omega_{j-1}}\h_{j}(\varphi_x +\psi)\overline{\psi}dx \right|&=&\displaystyle  \left|-\int_{\omega_{j-1}}(\h_{j}^\prime\overline{\psi}+\h_{j}\overline{\psi_x})\varphi dx+\int_{\omega_{j-1}}\h_{j}|\psi|^2 dx\right|\vspace{0.15cm}\\
\displaystyle &\leq& \displaystyle  \left(\max_{x\in \overline{\omega_j}} |\mathsf{h}_j^\prime(x)|\|\psi\|_{L^2(\omega_{j-1})}+\|\psi_x\|_{L^2(\omega_{j-1})}\right)\|\varphi\|_{L^2(\omega_{j-1})} + \|\psi\|^2_{L^2(\omega_{j-1})}. 
\end{array}
	$$
	Moreover, owing to Cauchy-Schwarz's inequality, we obtain
	$$
\left|	\int_{\omega_{j-1}}\h_{j}\theta_x \overline{\psi}dx\right| \leq  \|\theta_x\|_{L^2(\omega_{j-1})}\|\psi\|_{L^2(\omega_{j-1})}
	$$
	and
	$$
\left|	\int_{\omega_{j-1}}\h_{j}( \la^{-2}f^4+i \la^{-1}f^3)\overline{\psi}dx\right|\leq  \left(\la^{-2}+\la^{-1}\right)\|F\|_{\HH}\|\psi\|_{L^2(\omega_{j-1)}}.
	$$
	Using \eqref{4.25} and the above inequalities,  we deduce \eqref{4.21p}.\\\linebreak
\textbf{Step 3.} In this step, we prove that 
\begin{equation}\label{4.24}
\|\la \psi\|^2_{L^2(\omega_3)}=O\left(\la^{-2}\right) \ \ \text{and} \ \ \|\psi_{xx}\|^2_{L^2(\omega_3)}=O(1).
\end{equation}
First, taking $j=1$ in \eqref{4.20p}, using Lemma \ref{p3-1stlemps}, and \eqref{star}, we get 
\begin{equation}\label{4.26}
\|\psi_x\|^2_{L^2 (\omega_1)}=o(\la^{-1}).
\end{equation}
Now, take $j=2$ in \eqref{4.21p} and use the above estimation with Lemma \ref{p3-1stlemps} and \eqref{star}, we deduce that
\begin{equation}\label{4.27}
\|\la \psi\|^2_{L^2 (\omega_2)}=o(\la^{-1}).
\end{equation}
Multiplying \eqref{4.11p} and \eqref{4.12p} by $-ik_2 (\delta\la)^{-1}\h_{3}\overline{\psi_x}$ and $\h_{3}\overline{\psi}$, respectively, and integrating over $\omega_{2}$, we have
\begin{equation}\label{4.21}
k_2 \rho_3  \delta^{-1}\mathtt{I}_1 -ik_2 (\delta \la )^{-1}\mathtt{I}_2 +k_2 \int_{\omega_{2}}\h_{3}|\psi_x|^2 dx =-ik_2 \int_{\omega_{2}}\h_{3}( \delta^{-1}\rho_3 \la^{-3}f^5+i \la^{-3}f^3_x)\overline{\psi_x}dx 
\end{equation}
and
\begin{equation}\label{4.22}
\rho_2 \int_{\omega_{2}}\h_{3}|\la \psi|^2dx+k_2 \mathtt{I}_3-k_1  \int_{\omega_{2}}\h_{3}(\varphi_x+\psi)\overline{\psi}dx-\delta  \int_{\omega_{2}}\h_{3}\theta_x \overline{\psi}dx=-\rho_2  \int_{\omega_{2}}\h_{3} (\la^{-2}f^4+i \la^{-1}f^3)\overline{\psi}dx,
\end{equation}
where
$$
\mathtt{I}_1=\int_{\omega_{2}}\h_{3}\theta \overline{\psi_x}dx, \ \  \mathtt{I}_2= \int_{\omega_{2}}\h_{3}q_x \overline{\psi_x}dx  \ \text{and} \
\mathtt{I}_3 = \int_{\omega_{2}}\h_{3}\psi_{xx}\overline{\psi}dx.
$$
Using integration by parts, $\mathtt{I}_1$, $\mathtt{I}_2$, and $\mathtt{I}_3$ become
$$
\mathtt{I}_1 =-\int_{\omega_{2}}\h_{3}\theta_x \overline{\psi}dx - \int_{\omega_{2}}\h_{3}^\prime \theta \overline{\psi}dx, 
$$
$$
\mathtt{I}_2 = - \int_{\omega_{2}}\h_{3}q \overline{\psi_{xx}}dx- \int_{\omega_{2}}\h_{3}^\prime q \overline{\psi_x}dx
$$
and
\begin{equation}\label{I3}
\mathtt{I}_3 = -\int_{\omega_{2}}\h_{3}|\psi_{x}|^2 dx-\int_{\omega_{2}}\h_{3}^\prime \psi_x \overline{\psi}dx
\end{equation}
Now, Cauchy-Schwarz's inequality yields
$$
\begin{array}{lll}
\displaystyle|\mathtt{I}_1|&\leq& \displaystyle \underbrace{\|\sqrt{\h_{3}}\theta_x\|_{L^2(\omega_{2})}\|\sqrt{\h_{3}}\psi\|_{L^2(\omega_{2})}}_{=:A_1}+\max_{x\in \overline{\omega_3}} |\mathsf{h}_3^\prime(x)|\|\theta\|_{L^2(\omega_{2})}\|\psi\|_{L^2(\omega_{2})}.
\end{array}
$$
In view of  \eqref{4.11p}, we have
$$
k_2 \overline{\psi_{xx}}=-\rho_2 \la^2 \overline{\psi}+k_1 (\overline{\varphi_x}+\overline{\psi})+\delta \overline{\theta_x}-\rho_2 \la^{-2}\overline{f^4}+i\rho_2 \la^{-1}\overline{f^3} \ \ \text{in} \ \ \omega_{2}.
$$
Hence, 
$$
\begin{array}{lll}
\displaystyle k_2 \mathtt{I}_2& =& \displaystyle \rho_2 \la^2 \int_{\omega_{2}}\h_{3}q\overline{\psi}dx-k_1 \int_{\omega_{2}}\h_{3}q (\overline{\varphi_x}+\overline{\psi})dx-\delta \int_{2}\h_{3}q \overline{\theta_x}dx +\rho_2 \int_{\omega_{2}}\h_{3}q(\la^{-2}\overline{f^4}-i\la^{-1}\overline{f^3})dx \vspace{0.15cm}\\
&& \displaystyle -\, k_2\int_{\omega_{2}}\h^\prime_{3}q \overline{\psi_x}dx.
\end{array}
$$
Consequently, we obtain
$$
\begin{array}{lll}
\displaystyle k_2 \la^{-1} |\mathtt{I}_2|&\leq& \displaystyle \underbrace{\rho_2 \|\sqrt{\h_{3}}q\|_{L^2(\omega_{2})}\|\sqrt{\h_{3}}\la \psi\|_{L^2(\omega_{2})}}_{=:A_2}+k_1 \la^{-1}\|q\|_{L^2(\omega_{2})}\left(\|\varphi_x\|_{L^2(\omega_{2})}+\|\psi\|_{L^2(\omega_{2})}\right)\vspace{0.15cm}\\
&&+\, 
\displaystyle  \delta \la^{-1}\|q\|_{L^2(\omega_{2})}\|\theta_x\|_{L^2(\omega_{2})}+\rho_2 (\la^{-3}+\la^{-2})\|q\|_{L^2(\omega_{2})}\|F\|_{\HH}\vspace{0.15cm}\\ &&+\,  \displaystyle k_2\la^{-1}\max_{x\in \overline{\omega_3}} |\mathsf{h}_3^\prime(x)|\|q\|_{L^2(\omega_{2})}\|\psi_x\|_{L^2(\omega_{2})}.
\end{array}
$$
Adding \eqref{4.21} and \eqref{4.22}, then using \eqref{I3}, we infer
$$
\begin{array}{lll}
\rho_2 \|\sqrt{\h_{3}}\la \psi\|^2_{L^2(\omega_{2})}&=&\displaystyle\left|-\frac{k_2 \rho_3 }{\delta}\mathtt{I}_1 +\frac{ik_2}{ \delta \la}\mathtt{I}_2 +k_2 \int_{\omega_{2}}\h^\prime_{3}\psi_x \overline{\psi}dx + k_1  \int_{\omega_{2}}\h_{3}(\varphi_x+\psi)\overline{\psi}dx+\delta  \int_{\omega_{2}}\h_{3}\theta_x \overline{\psi}dx\right|\vspace{0.15cm}\\
\displaystyle &\leq& \displaystyle  k_2 \rho_3\delta^{-1}|\mathtt{I}_1|+\delta^{-1}k_2 \la^{-1} |\mathtt{I}_2|+k_2 \max_{x\in \overline{\omega_3}} |\mathsf{h}_3^\prime(x)| \|\psi_x \|_{L^2 (\omega_{2})}\|\psi \|_{L^2 (\omega_{2})}\vspace{0.15cm}\\
&&+\,  \displaystyle \underbrace{k_1 \|\sqrt{\mathsf{h}_{3}}\varphi_x \|_{L^2 (\omega_{2})}\|\sqrt{\mathsf{h}_{3}}\psi\|_{L^2 (\omega_{2})}}_{=:A_3}+k_1 \|\psi \|_{L^2 (\omega_{2})}^2+\underbrace{\delta \|\sqrt{\mathsf{h}_{3}}\theta_x \|_{L^2 (\omega_{2})}\|\sqrt{\mathsf{h}_{3}}\psi\|_{L^2 (\omega_{2})}}_{=:A_4}
\end{array}
$$
Thus, we get
$$
\begin{array}{lll}
\rho_2 \|\sqrt{\h_{3}}\la \psi\|^2_{L^2(\omega_{2})}&\leq & \displaystyle k_2 \rho_3 \delta^{-1}A_1 +k_2 \rho_3 \delta^{-1}\max_{x\in \overline{\omega_3}} |\mathsf{h}_3^\prime(x)|\|\theta\|_{L^2(\omega_{2})}\|\psi\|_{L^2(\omega_{2})}+\delta^{-1}A_2\vspace{0.15cm}\\&&+\,  \delta^{-1}k_1 \la^{-1}\|q\|_{L^2(\omega_{2})}\left(\|\varphi_x\|_{L^2(\omega_{2})}+\|\psi\|_{L^2(\omega_{2})}\right)+\la^{-1}\|q\|_{L^2(\omega_{2})}\|\theta_x\|_{L^2(\omega_{2})}\vspace{0.15cm}\\&&+\,  \displaystyle \delta^{-1}\rho_2 (\la^{-3}+\la^{-2})\|q\|_{L^2(\omega_{2})}\|F\|_{\HH}+\delta^{-1}k_2\la^{-1}\max_{x\in \overline{\omega_3}} |\mathsf{h}_3^\prime(x)|\|q\|_{L^2(\omega_{2})}\|\psi_x\|_{L^2(\omega_{2})}\vspace{0.15cm}\\
&&+\,  \displaystyle k_2 \max_{x\in \overline{\omega_3}} |\mathsf{h}_3^\prime(x)| \|\psi_x \|_{L^2 (\omega_{2})}\|\psi \|_{L^2 (\omega_{2})}+A_3+ k_1 \|\psi \|_{L^2 (\omega_{2})}^2+A_4,
\end{array}
$$
as a consequence of \eqref{star}, Lemma \ref{p3-1stlemps}, \eqref{4.26}, \eqref{4.27}, we deduce that
$$
\rho_2 \|\sqrt{\h_{3}}\la \psi\|^2_{L^2(\omega_{2})}\leq k_2 \rho_3 \delta^{-1}A_1 +\delta^{-1}A_2+A_3+A_4+o(\la^{-2}).
$$
Using Young's inequalities, we obtain 
$$
\left\{\begin{array}{lll}
\displaystyle A_1 \leq \displaystyle \frac{1}{2c_1(\la)}\|\sqrt{\h_{3}}\theta_x\|^2_{L^2(\omega_{2})}+\frac{c_1(\la)}{2}\|\sqrt{\h_{3}}\psi\|^2_{L^2(\omega_{2})},\vspace{0.15cm}
\\
A_2 \leq \displaystyle\frac{\rho_2}{2c_2} \|\sqrt{\h_{3}}q\|_{L^2(\omega_{2})}^2+\frac{c_2\rho_2}{2}\|\sqrt{\h_{3}}\la \psi\|_{L^2(\omega_{2})},\vspace{0.15cm}
\\
A_3\leq \displaystyle\frac{k_1^2}{2c_3(\la) \rho_2}\|\sqrt{\mathsf{h}_{3}}\varphi_x\|_{L^2 (\omega_{2})}^2+\frac{\rho_2 c_3 (\la)}{2}\|\sqrt{\mathsf{h}_{3}}\psi\|_{L^2 (\omega_{2})}^2,\vspace{0.15cm}
\\
A_4\leq \displaystyle\frac{\delta^2 }{2c_4(\la) \rho_2}\|\sqrt{\mathsf{h}_{3}}\theta_x\|_{L^2 (\omega_{2})}^2+\frac{\rho_2 c_4 (\la)}{2}\|\sqrt{\mathsf{h}_{3}}\psi\|_{L^2 (\omega_{2})}^2,
\end{array}\right.
$$
where $c_1(\la), c_2, c_3(\la), c_4(\la)$ are positive  constants to be determined. Hence we can infer that
$$
\begin{array}{lll}
\rho_2 \|\sqrt{\mathsf{h}_{3}}\la \psi\|^2_{L^2 (\omega_{2})}&\leq&\displaystyle  \frac{k_2 \rho_3 \delta^{-1}}{2c_1(\la)}\|\sqrt{\h_{3}}\theta_x\|^2_{L^2(\omega_{2})}+\frac{k_2 \rho_3 \delta^{-1}c_1(\la)}{2}\|\sqrt{\h_{3}}\psi\|^2_{L^2(\omega_{2})} +\frac{\delta^{-1}\rho_2}{2c_2} \|\sqrt{\h_{3}}q\|_{L^2(\omega_{2})}^2\vspace{0.15cm}\\
&&+\,  \displaystyle \frac{c_2 \rho_2 \delta^{-1}}{2}\|\sqrt{\h_{3}}\la \psi \|^2_{L^2(\omega_{2})}+\frac{k_1^2}{2c_3(\la) \rho_2}\|\sqrt{\mathsf{h}_{3}}\varphi_x\|_{L^2 (\omega_{2})}^2+\frac{\rho_2 c_3 (\la)}{2}\|\sqrt{\mathsf{h}_{3}}\psi\|_{L^2 (\omega_{2})}^2\vspace{0.15cm}\\
&&+\,  \displaystyle \frac{\delta^2 }{2c_4(\la) \rho_2}\|\sqrt{\mathsf{h}_{3}}\theta_x\|_{L^2 (\omega_{2})}^2+\frac{\rho_2 c_4 (\la)}{2}\|\sqrt{\mathsf{h}_{3}}\psi\|_{L^2 (\omega_{2})}^2+o(\la^{-2}).

\end{array}
$$
Now, taking 
$$
c_1(\la)=\frac{\la^2 \delta \rho_2}{4k_2 \rho_3 }, \ c_2=\frac{\delta}{4}, \ c_3(\la )=c_4(\la)=\frac{\la^2}{4},
$$
we get 
$$
\frac{\rho_2}{2} \|\sqrt{\mathsf{h}_{3}}\la \psi\|^2_{L^2 (\omega_{2})}\leq\frac{2}{\rho_2 \la^2} \left[\frac{(k_2 \rho_3)^2 }{\delta^2 }+\delta^2 \right]\|\theta_x\|^2_{L^2(\omega_{2})}+\frac{2\rho_2}{\delta^2}\|q\|^2_{L^2(\omega_2)}+ \frac{2k_1^2}{\la^2  \rho_2}\|\varphi_x\|_{L^2 (\omega_{2})}^2 +o(\la^{-2}).
$$
From the above estimation, Lemma \ref{p3-1stlemps}, \eqref{star}, and \eqref{star2}, we deduce that
$$
\|\la \psi\|^2_{L^2(\omega_3)}\leq  \|\sqrt{\mathsf{h}_{3}}\la \psi\|^2_{L^2 (\omega_{2})}=O(\la^{-2}) \ \text{and} \ \|\psi_{xx}\|^2_{L^2(\omega_3)}=O(1).
$$
\\
\textbf{Step 4:} In this step, we finish the proof, with this intention, take $j=4$ in \eqref{4.20p}, then using \eqref{4.24}, \eqref{star} and Lemma \ref{p3-1stlemps}, we obtain 
\begin{equation*}
\|\psi_x\|^2_{L^2(\omega_4)}=o(\la^{-2}).
\end{equation*}
Finally, by taking $j=5$ in \eqref{4.21p} and by using the above estimation with \eqref{star} and Lemma \ref{p3-1stlemps}, we get
$$
\|\la \psi\|^2_{L^2 (\omega_5)}=o(\la^{-2}).
$$
\end{proof}\\
\begin{lem}\label{4.3}
	\rm{The solution $\Phi=(\varphi,u,\psi,v,\theta,q,z)^{\top}\in D(\AA)$ of  \eqref{pm-f1p}-\eqref{pm-f7p} satisfies the following estimations }
\begin{equation*}\label{phi}
\|\varphi_x\|_{L^2 (\omega_6)}^2 =o(1) \ \ \text{and} \ \ \|\la \varphi\|^2_{L^2 (\omega_7)}=o(1).
\end{equation*}
\end{lem}
\begin{proof}
First, multiplying \eqref{4.11p} by $\h_6\overline{\varphi_x}$ and integrating over $\omega_5$, then using integration by parts, Cauchy-Schwarz's inequality, \eqref{star} and Lemmas \ref{p3-1stlemps}-\ref{4.2}, we deduce that
$$
\begin{array}{lll}
\displaystyle k_1 \int_{\omega_5}\h_6 |\varphi_x|^2 dx &=&\displaystyle \left|\int_{\omega_5}\h_6(\rho_2 \la^2 \psi+k_2 \psi_{xx}-k_1 \psi-\delta \theta_x+\rho_2 \la^{-2}f^4+i\rho_2 \la^{-1}f^3)\overline{\varphi_x}dx\right| \vspace{0.15cm}\\
&=&\displaystyle \left|\int_{\omega_5}\h_6(\rho_2 \la^2 \psi-k_1 \psi-\delta \theta_x+\rho_2 \la^{-2}f^4+i\rho_2 \la^{-1}f^3)\overline{\varphi_x}dx-k_2 \int_{\omega_5}(\h^\prime_6 \overline{\varphi_x}+\h_6 \overline{\varphi_{xx}})\psi_xdx\right|\vspace{0.15cm}\\
&\leq & \displaystyle \max \left( \rho_2,k_1,\delta\right)\left[\la^2 \|\psi\|_{L^2(\omega_5)}+\|\psi\|_{L^2(\omega_5)}+\|\theta_x\|_{L^2(\omega_5)}+(\la^{-2}+\la^{-1})\|F\|_{\HH}\right]\|\varphi_x\|_{L^2 (\omega_5 )}\vspace{0.15cm}\\
&& +\, \displaystyle k_2 \max\left( \max_{x\in \overline{\omega_6}} |\mathsf{h}_6^\prime(x)|, 1\right)\left(\|\varphi_x\|_{L^2(\omega_5)}+\|\varphi_{xx}\|_{L^2(\omega_5)}\right)\|\psi_{x}\|_{L^2(\omega_5)}\vspace{0.15cm}\\
&=&o(1).
\end{array}
$$
Thus, the first estimation in Lemma \ref{4.3} is proved . Now, multiplying \eqref{4.10p} by $\h_7 \overline{\varphi}$ and integrating over $\omega_6$, then using integration by parts, Cauchy-Schwarz's inequality, \eqref{star}, Lemmas \ref{p3-1stlemps}-\ref{4.2}, the first estimation of Lemma \ref{4.3}, we obtain
$$
\begin{array}{lll}
\displaystyle \rho_1 \int_{\omega_6}\h_7|\la \varphi|^2dx&=&\displaystyle \left|\int_{\omega_6}\h_7 (-k_1 (\varphi_x+\psi)_x -\rho_1 \la^{-2}f^2-i\rho_1 \la^{-1}f^1)\overline{\varphi}dx\right|\vspace{0.15cm}\\
&=& \displaystyle \left|\int_{\omega_6}\h_7 (-\rho_1 \la^{-2}f^2-i\rho_1 \la^{-1}f^1)\overline{\varphi}dx+k_1 \int_{\omega_6} (\h_7^\prime \overline{\varphi}+\h_7 \overline{\varphi_x})(\varphi_x+\psi)dx\right| \vspace{0.15cm}\\
&\leq&\displaystyle  \rho_1 (\la^{-2}+\la^{-1})\|F\|_{\HH}\|\varphi\|_{L^2 (\omega_6)}\vspace{0.15cm}\\
&&+\,  \displaystyle k_1 \max \left(\max_{x\in \overline{\omega_7}} |\mathsf{h}_7^\prime(x)|, 1\right)\left( \|\varphi\|_{L^2 (\omega_6)}+\|\varphi_x\|_{L^2 (\omega_6)}\right)\left(\|\varphi_x\|_{L^2 (\omega_6)}+\|\psi\|_{L^2 (\omega_6)}\right)\vspace{0.15cm}\\
&=& o(1).
\end{array}
$$
Finally, the above estimation finished the proof.
\end{proof}\\

\noindent \textbf{Proof of Theorem \ref{p5-polthm}:}
First, from Lemmas \ref{4.2}-\ref{4.3}, we deduce that
\begin{equation}\label{4.30}
\int_{\omega_7}\left(\rho_1 |\la \varphi|^2 +k_1 |\varphi_x|^2+\rho_2 |\la \psi|^2 +k_2 |\psi_x|^2\right)dx =o(1).
\end{equation}
	Now, let us fix a function $\chi \in C^1 ([0,\ell])$ such that $\chi(0)=\chi(\ell)=0$.
Multiplying \eqref{4.10p} and \eqref{4.11p} by $-2\chi \overline{\varphi_x}$ and $-2\chi \overline{\psi_x}$, respectively and integrating over $(0,\ell)$, then taking the real part, yield
$$
- \int_0^{\ell} \chi \left(\rho_1\left|\la \varphi\right|^2+k_1\left| \varphi_x\right|^2\right)_xdx -\Re\left[2k_1 \int_0^{\ell} \chi \psi_x \overline{\varphi_x}dx  \right]=\Re\left[ 2\int_0^{\ell}\chi\left(\rho_1 \la^{-2}f^2+i\rho_1 \la^{-1}f^1\right)\overline{\varphi_x}dx\right]
$$
and
$$
\begin{array}{lll}
	\displaystyle - \int_0^{\ell} \chi \left(\rho_2\left|\la \psi\right|^2+k_2\left| \psi_x\right|^2\right)_xdx +\Re\left[2k_1 \int_0^{\ell} \chi \varphi_x \overline{\psi_x}dx +2k_1 \int_0^{\ell} \chi \psi \overline{\psi_x}dx  \right]\vspace{0.15cm}\\
	\displaystyle +\, \Re\left[2\delta \int_0^{\ell_0}\chi \theta_x \overline{\psi_x}dx\right]=\Re\left[ 2\int_0^{\ell}\chi\left(\rho_2 \la^{-2}f^4+i\rho_2 \la^{-1}f^3\right)\overline{\psi_x}dx\right].
\end{array}
$$
Adding the above equations, then using integration by parts, we get 
$$
\begin{array}{ll}
	\displaystyle \int_0^{\ell} \chi^\prime \left(\rho_1\left|\la \varphi\right|^2+k_1\left| \varphi_x\right|^2+\rho_2\left|\la \psi\right|^2+k_2\left| \psi_x\right|^2\right)dx=\Re\left[-2k_1 \int_0^{\ell}\chi \psi\overline{\psi_x}dx-2\delta \int_0^{\ell_0} \chi \theta_x \overline{\psi_x}dx \right.\\\displaystyle 
	\left. +\,2\int_0^{\ell}\chi\left(\rho_1 \la^{-2}f^2+i\rho_1 \la^{-1}f^1\right)\overline{\varphi_x}dx+2\int_0^{\ell}\chi\left(\rho_2 \la^{-2}f^4+i\rho_2 \la^{-1}f^3\right)\overline{\psi_x}dx\right].
\end{array}
$$
Using Cauchy-Schwarz's inequality in the above equation, we deduce that
$$
\begin{array}{lll}
	\displaystyle\int_0^{\ell} \chi^\prime \left(\rho_1\left|\la \varphi\right|^2+k_1\left| \varphi_x\right|^2+\rho_2\left|\la \psi\right|^2+k_2\left| \psi_x\right|^2\right)dx&\leq& \displaystyle  \max_{x\in [0,\ell]}\chi(x)\left[2k_1 \|\psi\|_{L^2(0,\ell)}+2\delta \|\theta_x\|_{L^2(0,\ell_0)}\right.\vspace{0.15cm}\\
	\displaystyle &&+\,  \left.\displaystyle 2 (\rho_1 +\rho_2 )(\la^{-2}+\la^{-1})\|F\|_{\HH}\right]\|U\|_{\HH}.
\end{array}
$$
Now, by taking $\chi =x\chi_0+(x-\ell)\chi_1$ with $\chi_0,\chi_1 \in C^1 ([0,\ell],[0,1])$ such that 
$$
\left\{\begin{array}{lll}
\chi_0=1 \ \ \text{in}\ [0,7\varepsilon],\vspace{0.15cm}\\ \chi_0=0 \ \ \text{in} \ [\ell_0-7\varepsilon,\ell],\vspace{0.15cm} \\ \chi_1=0 \ \ \text{in} \ [0,7\varepsilon],\vspace{0.15cm} \\ \chi_1=1 \ \ \text{in} \ [\ell_0-7\varepsilon,\ell],
\end{array}\right.
$$
we deduce that $\chi^\prime=\chi_0+\chi_1+x\chi_0^\prime+(x-\ell)\chi_1^\prime$ and
$$
\begin{array}{l}
\displaystyle\int_{(0,\ell)\backslash \omega_7 } \left(\rho_1\left|\la \varphi\right|^2+k_1\left| \varphi_x\right|^2+\rho_2\left|\la \psi\right|^2+k_2\left| \psi_x\right|^2\right)dx\vspace{0.15cm}\\
\leq \displaystyle \left[ (\ell_0-7\varepsilon)\max_{x\in \overline{ \omega_7}}|\chi_0^\prime (x)|+(\ell-7\varepsilon)\max_{x\in \overline{ \omega_7}}|\chi_1^\prime (x)|\right]\int_{\omega_7 } \left(\rho_1\left|\la \varphi\right|^2+k_1\left| \varphi_x\right|^2+\rho_2\left|\la \psi\right|^2+k_2\left| \psi_x\right|^2\right)dx\vspace{0.15cm}\\
\quad +\displaystyle  \max_{x\in [0,\ell]}\chi(x)\left[2k_1 \|\psi\|_{L^2(0,\ell)}+2\delta \|\theta_x\|_{L^2(0,\ell_0)}
\displaystyle +2 (\rho_1 +\rho_2 )(\la^{-2}+\la^{-1})\|F\|_{\HH}\right]\|\Phi\|_{\HH},
\end{array}
$$
consequently, from Lemmas \ref{p3-1stlemps}, \eqref{star}, and \eqref{4.30}, we get 
\begin{equation}\label{4.31}
\int_{(0,\ell)\backslash \omega_7 } \left(\rho_1\left|\la \varphi\right|^2+k_1\left| \varphi_x\right|^2+\rho_2\left|\la \psi\right|^2+k_2\left| \psi_x\right|^2\right)dx=o(1).
\end{equation}
Now, from \eqref{4.30} and \eqref{4.31}, we obtain
\begin{equation}\label{4.32}
	\int_0^\ell \left(\rho_1\left|\la \varphi\right|^2+k_1\left| \varphi_x\right|^2+\rho_2\left|\la \psi\right|^2+k_2\left| \psi_x\right|^2\right)dx=o(1).
\end{equation}
Finally, from the above estimation and Lemma \ref{p3-1stlemps}, we deduce that 
$$
\|\Phi\|_{\HH}=o(1).
$$
\xqed{$\square$}

\appendix
\section{Some notions and stability theorems}\label{p2-appendix}
In order to make this paper more self-contained, we recall in this short appendix some notions and stability results used in this work. 
\begin{defi}\label{App-Definition-A.1}{\rm
		Assume that $A$ is the generator of $C_0-$semigroup of contractions $\left(e^{tA}\right)_{t\geq0}$ on a Hilbert space $H$. The $C_0-$semigroup $\left(e^{tA}\right)_{t\geq0}$ is said to be 
		\begin{enumerate}
			\item[$(1)$] Strongly stable if 
			$$
			\lim_{t\to +\infty} \|e^{tA}x_0\|_H=0,\quad \forall\, x_0\in H.
			$$
			\item[$(2)$] Exponentially (or uniformly) stable if there exists two positive constants $M$ and $\varepsilon$ such that 
			$$
			\|e^{tA}x_0\|_{H}\leq Me^{-\varepsilon t}\|x_0\|_{H},\quad \forall\, t>0,\ \forall\, x_0\in H.
			$$
			\item[$(3)$] Polynomially stable if there exists two positive constants $C$ and $\alpha$ such that 
			$$
			\|e^{tA}x_0\|_{H}\leq Ct^{-\alpha}\|A x_0\|_{H},\quad \forall\, t>0,\ \forall\, x_0\in D(A).
			$$
			\xqed{$\square$}
	\end{enumerate}}
\end{defi}
\noindent To show  the strong stability of a $C_0$-semigroup  we rely on the following result due to Arendt-Batty \cite{Arendt01}. 
\begin{Theorem}\label{App-Theorem-A.2}{\rm
		{Assume that $A$ is the generator of a C$_0-$semigroup of contractions $\left(e^{tA}\right)_{t\geq0}$  on a Hilbert space $H$. If $A$ has no pure imaginary eigenvalues and  $\sigma\left(A\right)\cap i\mathbb{R}$ is countable,
			where $\sigma\left(A\right)$ denotes the spectrum of $A$, then the $C_0$-semigroup $\left(e^{tA}\right)_{t\geq0}$  is strongly stable.}\xqed{$\square$}}
\end{Theorem}
\noindent Finally for the  polynomial stability   of a $C_0-$semigroup of contractions we use the following result due to Borichev and Tomilov \cite{Borichev01} (see also \cite{Batty01} and \cite{RaoLiu01}).
\begin{Theorem}\label{bt}
	{\rm
		Assume that $A$ is the generator of a strongly continuous semigroup of contractions $\left(e^{tA}\right)_{t\geq0}$  on $H$.   If   $ i\mathbb{R}\subset \rho(A)$, then for a fixed $\ell>0$ the following conditions are equivalent
		\begin{equation}\label{h1}
			\limsup_{\la\in \R,\ \abs{\la}\rightarrow \infty}\frac{1}{|\la|^{\ell}}\|(i\la I-A)^{-1}\|_{\mathcal{L}(H)}<\infty,
		\end{equation}
		\begin{equation}\label{h2}
			\|e^{tA}\Phi_{0}\|^2_{H} \leq \frac{C}{t^{\frac{2}{\ell}}}\|\Phi_0\|^2_{D(A)},\hspace{0.1cm}\forall t>0,\hspace{0.1cm} \Phi_0\in D(A),\hspace{0.1cm} \text{for some}\hspace{0.1cm} C>0.
		\end{equation}\xqed{$\square$}}
\end{Theorem}


\end{document}